\documentclass{amsart}

\newtheorem{theorem}{Theorem}[section]
\newtheorem{lemma}[theorem]{Lemma}
\newtheorem{proposition}[theorem]{Proposition}
\newtheorem{corollary}[theorem]{Corollary}
\newtheorem{question}[theorem]{Question}
\newtheorem{problem}[theorem]{Problem}

\theoremstyle{definition}

\newtheorem{example}[theorem]{Example}
\newtheorem{examples}[theorem]{Examples}
\newtheorem{remark}[theorem]{Remark}
\newtheorem{remarks}[theorem]{Remarks}
\newtheorem{algorithm}[theorem]{Algorithm}

\usepackage{amscd,amssymb}

\begin{document}

\title[Weitzenb\"ock derivations and unipotent transformations]
{Constants of Weitzenb\"ock derivations\\
and invariants of unipotent transformations\\
acting on relatively free algebras}

\author[Vesselin Drensky and C. K. Gupta]
{Vesselin Drensky and C. K. Gupta}
\address{Institute of Mathematics and Informatics,
Bulgarian Academy of Sciences,
          1113 Sofia, Bulgaria}
\email{drensky@math.bas.bg}
\address{Department of Mathematics, University of Manitoba,
          Winnipeg, Manitoba, R3T 2N2 Canada}
\email{cgupta@cc.umanitoba.ca}

\thanks
{The work of the first author was partially supported by Grant MM-1106/2001 of
the Bulgarian Foundation for Scientific Research.}

\thanks
{The work of the second author was supported by NSERC, Canada.}


\begin{abstract} In commutative algebra, a Weitzenb\"ock derivation is
a nonzero triangular linear derivation of the polynomial algebra
$K[x_1,\ldots,x_m]$ in several variables over a field $K$ of characteristic 0.
The classical theorem of Weitzenb\"ock states that the algebra of constants
is finitely generated. (This algebra coincides with
the algebra of invariants of a single unipotent transformation.)
In this paper we study the problem of finite generation of the algebras
of constants of triangular linear derivations of finitely generated
(not necessarily commutative or associative) algebras over $K$
assuming that the algebras are free in some sense
(in most of the cases relatively free algebras in varieties of
associative or Lie algebras). In this case the algebra of constants also
coincides with the algebra of invariants of some unipotent transformation.
\par
The main results are the following: 1. We show that the subalgebra of constants
of a factor algebra can be lifted to the subalgebra of constants.
2. For all varieties of associative algebras which are not nilpotent in Lie sense
the subalgebras of constants of the relatively free algebras of rank $\geq 2$
are not finitely generated. 3.
We describe the generators of the subalgebra of constants
for all factor algebras $K\langle x,y\rangle/I$
modulo a $GL_2(K)$-invariant ideal $I$.
4. Applying known results from commutative algebra, we construct classes
of automorphisms of the algebra generated by two generic $2\times 2$ matrices.
We obtain also some partial results on relatively free Lie algebras.
\end{abstract}

\maketitle

\section{Introduction}

We fix a base field $K$ of characteristic 0, an integer $m\geq 2$
and a set of symbols $X=\{x_1,\ldots,x_m\}$. We call the elements of $X$
variables. Sometimes we shall use other symbols, e.g. $y,z,y_i$, etc.
for the elements of $X$. We denote by $V_m$ the vector space with
basis $X$.
\par
Let $K\langle X\rangle=K\langle x_1,\ldots,x_m\rangle$ be the free
unitary associative algebra freely generated by $X$ over $K$. The elements
of $K\langle X\rangle$ are linear combinations of words $x_{j_1}\cdots x_{j_n}$
in the noncommuting variables $X$. The general linear group $GL_m=GL_m(K)$
acts naturally on the vector space $V_m$ and this action is extended diagonally
on $K\langle X\rangle$ by the rule
\[
g(x_{j_1}\cdots x_{j_n})=g(x_{j_1})\cdots g(x_{j_n}),\
g\in GL_m,\ x_{j_1},\ldots,x_{j_n}\in X.
\]
All associative algebras which we consider in this paper are homomorphic
images of $K\langle X\rangle$ modulo ideals $I$ which are invariant under
this $GL_m$-action. We shall use the same symbols $x_j$ and $X$ for the generators
and the whole generating set of $K\langle X\rangle/I$.
Most of the algebras in our considerations will
be relatively free algebras in varieties of unitary associative algebras.
Examples of relatively free algebras are the polynomial algebra
$K[X]$ and the free algebra $K\langle X\rangle$ which are free, respectively,
in the varieties of all commutative algebras and all associative algebras.
We also shall consider Lie algebras which are
homomorphic images of the free Lie algebra with $X$ as a free generating set
modulo ideals which are also $GL_m$-invariant.
\par

Let $A$ be any (not necessarily associative or Lie) algebra over $K$.
Recall that the $K$-linear operator $\delta$ acting on $A$ is called a
derivation of $A$ if
\[
\delta(uv)=\delta(u)v+u\delta(v)\ \text{\rm for all }\ u,v\in A.
\]
The elements $u\in A$ which belong to the kernel of $\delta$ are called
constants of $\delta$ and form a subalgebra of $A$ which we shall denote by
$A^{\delta}$. The derivation $\delta$ is locally nilpotent if for any
$u\in A$ there exists a positive integer $n$ such that $\delta^n(u)=0$.
If $\delta$ is a locally nilpotent derivation of $A$, then the linear operator
of $A$
\[
\exp\delta=1+\frac{\delta}{1!}+\frac{\delta^2}{2!}+\cdots
\]
is well defined and is an automorphism of the $K$-algebra $A$.
It is easy to see that $A^{\delta}$ coincides with the subalgebra of fixed
points (or invariants) of $\exp\delta$ which we shall denote by
$A^{\exp\delta}$. The mapping $\alpha\to\exp(\alpha\delta)$,
$\alpha\in K$, defines an additive action of $K$ on $A$.
It is well known that for polynomial algebras every additive action of $K$
is of this kind, see for more details Snow \cite{Sn}.
See also Drensky and Yu \cite{DY1} for
relations between exponents of locally nilpotent derivations
and automorphisms $\varphi$ with the property that
the orbit $\{\varphi^n(a)\mid n\in{\mathbb Z}\}$ of each $a\in A$
spans a finite dimensional vector space in
the noncommutative case.
\par

If $A=K\langle X\rangle/I$ for some $GL_m$-invariant ideal $I$, then
the derivation $\delta$ of $A$ is called triangular, if $\delta(x_j)$,
$j=1,\ldots,m$, belongs to the subalgebra of $A$ generated by $x_1,\ldots,x_{j-1}$.
Clearly, the triangular derivations are locally nilpotent.
If $\delta$ acts linearly on the vector space $V_m=\sum_{j=1}^mKx_j\subset A$,
then it is called linear.
\par
If $\delta$ is a triangular derivation, then $\exp\delta$ is a triangular
automorphism of $A$, with the property
\[
\exp\delta(x_j)=x_j+f_j(x_1,\ldots,x_{j-1}),\ j=1,\ldots,m,
\]
where $f_j(x_1,\ldots,x_{j-1})$ depends on $x_1,\ldots,x_{j-1}$ only.
Every triangular automorphism $\varphi$ of this form can be presented
in the form $\varphi=\exp\delta$ for some triangular derivation
\[
\delta=\log(\varphi)=\frac{\varphi-1}{1}-\frac{(\varphi-1)^2}{2}
+\frac{(\varphi-1)^3}{3}-\frac{(\varphi-1)^4}{4}+\cdots.
\]
(The $K$-linear operator $\delta$ of $A$ is well defined because
the linear operators $(\varphi-1)^k$ map
every $f(x_1,\ldots,x_j)\in A$ to a polynomial
depending on $x_1,\ldots,x_j$ only,
$\text{deg}_{x_j}(\varphi-1)^kf\leq \text{deg}_{x_j}f-k$
and $(\varphi-1)K=0$.)
\par
Every locally nilpotent linear derivation $\delta$
is triangular with respect to a suitable basis of $V_m$ and the
automorphism $\exp\delta$ is a unipotent linear transformation (i.e.
an automorphism of the algebra $A$ which acts as a unipotent linear
operator on $V_m$).
\par

In commutative algebra, the triangular linear derivations
of the polynomial algebra $K[X]=K[x_1,\ldots,x_m]$ are called
Weitzenb\"ock derivations.
The classical theorem of Weitzenb\"ock \cite{W} states that
the algebra of constants of such a derivation
is finitely generated. This algebra coincides with
the algebra of invariants of a single unipotent transformation.
\par
In this paper we study the problem of finite generation of the algebras
of constants of triangular linear derivations of (usually noncommutative)
algebras $K\langle X\rangle/I$ where the ideal $I$ is $GL_m$-invariant.
As in the commutative case, the algebra of constants
coincides with the algebra of invariants of some unipotent transformation.
The paper is organized as follows. Below we assume that $\delta$
is a nonzero triangular linear derivation of $K\langle X\rangle$ which induces
a derivation (which we shall also denote by $\delta$)
on the factor algebras of $K\langle X\rangle$ modulo $GL_m$-invariant ideals.
\par
In Section 2 we present a short survey on constants of locally nilpotent derivations
and invariant theory both in the commutative and noncommutative case,
giving some motivation for our investigations. We believe that some of the
results exposed there can serve as a motivation and inspiration for further
investigations on noncommutative algebras.
\par
Section 3 presents a summary of the results on the Weitzenb\"ock derivations
of polynomial algebras which we need in the next sections.
\par
In Section 4 we are interested in the problem of lifting the constants:
If $I\subset J$ are two $GL_m$-invariant ideals of $K\langle X\rangle$, then we
show that the subalgebra of constants $(K\langle X\rangle/J)^{\delta}$
can be lifted to the subalgebra of constants $(K\langle X\rangle/I)^{\delta}$.
In the special case of algebras with two generators $x,y$
we may assume that $\delta(x)=0$, $\delta(y)=x$. Then
the subalgebra of constants is spanned by elements which have a very special
behaviour under the action of the general linear group $GL_2$, the so called
highest weight vectors. This allows to involve classical combinatorial techniques
as theory of generating functions and representation theory of general linear
groups.
\par
In Section 5 we present various examples of subalgebras of constants
of relatively free associative algebras. In particular, he handle the case
of the free algebra $K\langle x,y\rangle$ and show that the algebra
of constants is generated by $x$ and a set of $SL_2(K)$-invariants
which we describe explicitly. As a consequence, we obtain a similar
generating set for all factor algebras $K\langle x,y\rangle/I$.
\par
Section 6 considers relatively free algebras $F_m({\mathfrak W})$
in varieties $\mathfrak W$ of
associative algebras. It is known that every variety $\mathfrak W$
is either nilpotent in Lie sense or contains the algebra of $2\times 2$
upper triangular matrices. We show that
for all $\mathfrak W$ which are not nilpotent in Lie sense
the subalgebras of constants $F_m({\mathfrak W})^{\delta}$
are not finitely generated.
\par
In Section 7 we apply results from commutative algebra and construct classes
of automorphisms of the relatively free algebra
$F_2(\text{\rm var }M_2(K))$. This algebra is isomorphic to the
algebra generated by two generic $2\times 2$ matrices $x$ and $y$.
The centre of the associated generic trace algebra
(which coincides with the algebra of invariants of two $2\times 2$ matrices
under simultaneous conjugation by $GL_2$) is generated by the traces
of $x$, $y$ and $xy$ and the determinants of $x$ and $y$ and is
isomorphic to the polynomial algebra in five variables.
We want to mention that up till now
most of the investigations have been performed in the opposite direction.
The automorphisms of
$F_2(\text{\rm var }M_2(K))$ and of the trace algebra
have been used to produce automorphisms of the polynomial algebra
in five variables, see e.g. Bergman \cite{B}, Alev and Le Bruyn \cite{AL},
Drensky and Gupta \cite{DG}.
\par
Finally, we obtain also some partial results on relatively free Lie algebras.

\section{Survey}

\subsection{Motivation from Commutative Algebra}
Locally nilpotent derivations of the polynomial algebra
$K[X]=K[x_1,\ldots,x_m]$ have been studied for many decades and have had
significant impact on different branches of algebra and invariant theory,
see e.g. the books by Nowicki \cite{No} and van den Essen \cite{E2}.
\par
Let $G$ be a subgroup of $GL_m$ and let $K[X]^G=K[x_1,\ldots,x_m]^G$
be the algebra of $G$-invariants. The problem for finite generation
of $K[X]^G$ was the main motivation for the famous Hilbert Fourteenth Problem
\cite{H}. The theorem of Emmy Noether \cite{N} gives the finite generation
of $K[X]^G$ for finite groups $G$. More generally, the Hilbert-Nagata theorem
states the finite generation of $K[X]^G$ for reductive groups $G$, see e.g. \cite{DC}.
\par
The first counterexample of Nagata \cite{N1}
to the Hilbert Fourteenth Problem was the non-finitely generated
algebra of invariants
$K[x_1,\ldots,x_{32}]^G$ of a specially constructed
triangular linear group $G$. Today, most of the known counterexamples have been
obtained (or can be obtained) as algebras of constants of some
derivations. This includes the original counterexample of Nagata,
see Derksen \cite{De} who was the first to recognize the connection between
the Hilbert 14-th problem and constants of derivations (but his derivations
were not always locally nilpotent) and the counterexample of Daigle and
Freudenburg \cite{DF} of a triangular (but not linear) derivation
of $K[x_1,\ldots,x_5]$ with not finitely generated algebra of constants.
For more counterexamples to the Hilbert 14-th problem
we refer to the recent survey by Freudenburg \cite{Fr2}.
\par
The theorem of Weitzenb\"ock gives the finite generation of the algebra
of constants for a triangular linear derivation or, equivalently, for
the algebra of invariants of a single unipotent transformation.
(This contrasts to the counterexample of Nagata described above.)
The original proof of Weitzenb\"ock from 1932 was for $K={\mathbb C}$.
Later  Seshadri \cite{Se} found a proof for any field $K$
of charactersitic 0. A simple proof for $K=\mathbb C$ using ideas from
\cite{Se} has been recently given by Tyc \cite{T}. To the best of our knowledge,
no constructive proof, with effective estimates of the degree of the generators
of the algebra of constants has been given up till now.
\par
For each dimension $m$ there are only finite number of
essentially different Weitzen\-b\"ock derivations to study:
Up to a linear change of the coordinates, the Weitzen\-b\"ock derivations $\delta$
are in one-to-one correspondence with the partition
$(p_1+1,p_2+1,\ldots,p_s+1)$ of $m$, where $p_1\geq p_2\geq\cdots\geq p_s\geq 0$,
$(p_1+1)+(p_2+1)+\cdots+(p_s+1)=m$, and the correspondence is given
in terms of the Jordan normal form $J(\delta)$ of the matrix of the derivation
\[
J(\delta)=\begin{pmatrix}
J_1&0&\cdots&0\\
0&J_2&\cdots&0\\
\vdots&\vdots&\cdots&\vdots\\
0&0&\cdots&J_s
\end{pmatrix},\quad\ \text{ \rm where }\
J_i=\begin{pmatrix}
0&1&0&\cdots&0&0\\
0&0&1&\cdots&0&0\\
\vdots&\vdots&\vdots&\cdots&\vdots&\vdots\\
0&0&0&\cdots&0&1\\
0&0&0&\cdots&0&0
\end{pmatrix},
\]
is the $(p_i+1)\times(p_i+1)$ Jordan cell with zero diagonal.
\par

Another important application of locally nilpotent derivations
is the construction of candidates for wild automorphisms of polynomial
algebras, see e.g. the survey of Drensky and Yu \cite{DY2}.
Typical example is the following. If $\delta$ is a Weitzenb\"ock derivation
of $K[x_1,\ldots,x_m]$ and $0\not= w\in K[x_1,\ldots,x_m]^{\delta}$, then
$\Delta=w\delta$ is also a locally nilpotent derivation of
$K[x_1,\ldots,x_m]$ with the same algebra of constants as $\delta$ and
$\exp\Delta$ is an automorphism of $K[x_1,\ldots,x_m]$.
By the theorem of Martha Smith \cite{Sm}, all such automorphisms are stably
tame and become tame if extended to $K[x_1,\ldots,x_m,x_{m+1}]$ by
$(\exp\Delta)(x_{m+1})=x_{m+1}$. The famous Nagata automorphism
of $K[x,y,z]$, see \cite{N2},
also can be obtained in this way:
We define the derivation $\delta$ by
\[
\delta(x)=-2y,\quad \delta(y)=z,\quad \delta(z)=0,\quad
w=xz+y^2\in K[x,y,z]^{\delta},
\]
and for $\Delta=w\delta$ the Nagata automorphism is $\nu=\exp\Delta$:
\begin{align*}
\nu(x)&=x+(-2y)\frac{w}{1!}+(-2z)\frac{w^2}{2!}
=x-2(xz+y^2)y-(xz+y^2)^2z,\\
\nu(y)&=y+z\frac{w}{1!}
=y+(xz+y^2)z,\\
\nu(z)&=z.
\end{align*}
Recently Shestakov and Umirbaev \cite{SU} proved that the Nagata automorphism
is wild. It is interesting to mention that their approach is based on
Poisson algebras and methods of noncommutative, and even nonassociative, algebras.
\par
There are few exceptions of locally nilpotent derivations and their exponents
which do not arrise immediately from triangular derivations: the derivations of
Freudenburg (obtained with his local slice construction \cite{Fr1})
and the automorphisms of Drensky and Gupta
(obtained by methods of noncommutative algebra, \cite{DG}). Later,
Drensky, van den Essen and Stefanov \cite{DES} have shown that the
automorphisms from \cite{DG}
also can be obtained in terms of locally nilpotent derivations
and are stably tame.

\subsection{Noncommutative Invariant Theory}
An important part of noncommutative invariant theory is devoted to
the study of the algebra of invariants of a linear group
$G\subset GL_m$ acting on the free associative algebra
$K\langle X\rangle=K\langle x_1,\ldots,x_m\rangle$,
relatively free algebras $F_m({\mathfrak W})$ in varieties of
associative algebras $\mathfrak W$, the free Lie algebra $L_m=L(X)$
and relatively free algebras $L_m({\mathfrak V})$ in varieties
of Lie algebras $\mathfrak V$. For more detailed exposition we refer to the
surveys on noncommutative invariant theory
by Formanek \cite{F1}, Drensky \cite{D5} and the survey on algorithmic methods
for relatively free semigroups, groups and algebras by
Kharlampovich and Sapir \cite{KS}.

\subsubsection{Free Associative Algebras}
By a theorem of Lane \cite{L} and Kharchenko \cite{K1},
the algebra of invariants $K\langle X\rangle^G$ is always a free algebra
(independently of the properties of $G\subset GL_m$).
By the theorem of Dicks and Formanek \cite{DiF} and Kharchenko \cite{K1},
if $G$ is finite, then $K\langle X\rangle^G$ is finitely generated if and only if
$G$ is cyclic and acts on $V_m=\sum_{j=1}^mKx_j$ as a group of scalar multiplications.
This result was generalized for a much larger class of groups by Koryukin \cite{Ko}
who also established a finite generation of $K\langle X\rangle^G$ if we equip it
with a proper action of the symmetric group.
\par
Recall that if $V$ is a multigraded vector space which is a direct sum
of its multihomogeneous components $V^{(n_1,\ldots,n_m)}$, then the Hilbert
series of $V$ is defined as the formal power series
\[
H(V,t_1,\ldots,t_m)=\sum\text{\rm dim}(V^{(n_1,\ldots,n_m)})t_1^{n_1}\cdots t_m^{n_m}.
\]
If $V$ is ``only'' graded with homogeneous components $V^{(n)}$,
then its Hilbert series is
\[
H(V,t)=\sum_{n\geq 0}\text{\rm dim}(V^{(n)})t^n.
\]
Dicks and Formanek \cite{DiF} proved also an analogue
of the Molien formula for the Hilbert series of $K\langle X\rangle^G$,
$\vert G\vert<\infty$, which was generalized for compact groups $G$ by
Almkvist, Dicks and Formanek \cite{ADF} (an analogue of the Molien-Weyl formula
in classical invariant theory). In particular, Almkvist, Dicks and Formanek
showed that the Hilbert series of the algebra of invariants
$K\langle X\rangle^g$ is an algebraic function if $g$ is a unipotent matrix.
(Hence the same holds for the algebra of constants $K\langle X\rangle^{\delta}$
for a Weitzenb\"ock derivation $\delta$.)

\subsubsection{Relatively Free Associative Algebras}
Let $f(x_1,\ldots,x_m)\in K\langle x_1,x_2,\ldots\rangle$ be an element of
the free algebra of countable rank. Recall that $f(x_1,\ldots,x_m)=0$
is a polynomial identity for the algebra $A$ if
$f(a_1,\ldots,a_m)=0$ for all $a_1,\ldots,a_m\in A$. The algebra is called PI,
if it satisfies some nontrivial polynomial identity.
The class of all algebras satisfying a given set
$U\subset K\langle x_1,x_2,\ldots\rangle$
of polynomial identities is called the variety of associative algebras
defined by the system $U$. We shall denote the varieties by German letters.
If $\mathfrak W$ is a variety, then $T({\mathfrak W})$ is the ideal
of $K\langle x_1,x_2,\ldots\rangle$ consisting of all polynomial identities of $\mathfrak W$
and the algebra
\[
F_m({\mathfrak W})=K\langle x_1,\ldots,x_m\rangle/
(K\langle x_1,\ldots,x_m\rangle\cap T({\mathfrak W}))
\]
is the relatively free algebra of rank $m$ in $\mathfrak W$.
The ideals $K\langle x_1,\ldots,x_m\rangle\cap T({\mathfrak W})$ of
$K\langle x_1,\ldots,x_m\rangle$ are invariant under all endomorphisms
of $K\langle x_1,\ldots,x_m\rangle$ and, in particular, are $GL_m$-invariant.
\par
Most of the work on invariant theory of relatively free algebras is devoted to the
description of the varieties $\mathfrak W$ such that $F_m({\mathfrak W})^G$
is finitely generated for all $m=2,3,\ldots$, and all groups
$G\subset GL_m$ from a given class $\mathfrak G$.
The description of such varieties for the
class of all finite groups is given in different terms by several authors,
see the surveys by Formanek \cite{F1}, Drensky \cite{D5}, Kharlampovich and Sapir
\cite{KS}. In particular, the finite generation of $F_m({\mathfrak W})^G$ for
all finite groups holds if and only if all finitely generated algebras of
$\mathfrak W$ are weakly noetherian (i.e. noetherian with respect to two-sided
ideals) which is equivalent to the fact that $\mathfrak W$ satisfies a polynomial
identity of a very special form. One of the simplest descriptions is
the following (see \cite{D3}): $F_m({\mathfrak W})^G$ is finitely generated for
all $m\geq 2$ and
all finite groups $G\subset GL_m$
if and only if $F_2({\mathfrak W})^g$ is finitely generated
for the linear transformation $g$ defined by $g(x_1)=-x_1$, $g(x_2)=x_2$.
\par
If we consider the finite generation of $F_m({\mathfrak W})^G$ for the class
all reductive groups $G$, then the results of Vonessen \cite{V},
Domokos and Drensky \cite{DD} give that $F_m({\mathfrak W})^G$ is finitely
generated for all reductive $G$ if and only if the finitely generated algebras
in $\mathfrak W$ are one-side noetherian. For unitary algebras this means that
$\mathfrak W$ satisfies the Engel identity $[x_2,x_1,\ldots,x_1]=0$.
\par
Concerning the Hilbert series of subalgebras of invariants of
relatively free algebras, Formanek \cite{F1} generalized the Molien-Weyl
formula for the Hilbert series of $K[x_1,\ldots,x_m]^G$ for $G$ finite or
compact to the case of any relatively free algebra, expressing
the Hilbert series of $F_m({\mathfrak W})^G$ in terms of the Hilbert series
$H(F_m({\mathfrak W}),t_1,\ldots,t_m)$. If $G$ is finite, then
$H(F_m({\mathfrak W})^G,t)$ involves the eigenvalues of all $g\in G$.
By a theorem of Belov \cite{Be}, the Hilbert series of
$F_m({\mathfrak W})$ is always a rational function and this imlies that
$H(F_m({\mathfrak W})^G,t)$ is also rational for $G$ finite.
For reductive $G$ the rationality of $H(F_m({\mathfrak W})^G,t)$
is known only for varieties $\mathfrak W$ satisfying a nonmatrix
polynomial identity, see Domokos and Drensky \cite{DD}.

\subsubsection{Lie Algebras}
We shall mention few results only. By a theorem of Bryant \cite{Br},
if $G$ is a nontrivial finite linear group, then the algebra of fixed points
of the free Lie algebra $L_m^G$ is never finitely generated.
This result was extended by Densky \cite{D4} to the fixed points
of all relatively free algebras $L_m({\mathfrak V})$ (and also for all finite
$G\not=1$) for nonnilpotent varieties $\mathfrak V$ of Lie algebras.
We refer also to the work done by several authors and mainly by
Bryant, Kovacz and St\"ohr about fixed points of free Lie algebras
in the modular case, see e.g. \cite{BKS} and the references there.

\subsection{Derivations of Free Algebras}
The algebra of constants of the formal partial derivatives
$\partial/\partial x_j$, $j=1,\ldots,m$, for
$K\langle X\rangle=K\langle x_1,\ldots,x_m\rangle$ was described by
Falk \cite{Fa}.
It is generated by all Lie commutators
$[[\ldots[x_{j_1},x_{j_2}],\ldots],x_{j_n}]$, $n\geq 2$. Specht \cite{Sp}
applied products of such commutators in the study of algebras with polynomial
identities, see also Drensky \cite{D2} or the book \cite{D6} for further
application to the theory of PI-algebras. It is known, see Gerritzen \cite{G},
that in this case the algebra of constants is free,
see also Drensky and Kasparian \cite{DK} for an explicit basis.
(The freedom of the algebra
of constants of the partial derivatives of $K[X]$ does not
follow immediately  from the result
of Lane \cite{L} and Kharchenko \cite{K1}. The derivations
$\partial/\partial x_j$ are locally nilpotent and their exponents
$\exp(\partial/\partial x_j)$ generate a group of automorphisms of
$K\langle X\rangle$ which consists of all translations of the form $x_i\to x_i+a_i$,
$a_i\in\mathbb Z$.
Although this group is a subgroup of the affine group,
we cannot apply directly \cite{L} and \cite{K1}
because the group is not linear.)
\par
Similar study of the algebra of constants in a very large class of
(not only associative) algebras was performed by Gerritzen and Holtkamp
\cite{GH} and Drensky and Holtkamp \cite{DH}.
We shall finish the survey section with the following, probably folklore known lemma.

\begin{lemma} \label{extending derivations}
Let $\mathfrak W$ be any variety of algebras and let
$F({\mathfrak W})$ be the relatively free algebra of any rank.
Every mapping from the free generating set to $F({\mathfrak W})$
can be extended to a derivation.
\end{lemma}

\begin{proof} We shall prove the lemma for relatively free associative algebras
of finite or countable rank only.
The same considerations work in the case of any infinite rank.
Let $\delta_0:\{x_1,x_2,\ldots\}\to F_{\infty}({\mathfrak W})$ be any mapping
and let $T({\mathfrak W})$ be the T-ideal of $K\langle x_1,x_2,\ldots\rangle$
of all polynomial identities of $\mathfrak W$. We fix
$f_1,\ldots,f_m\in K\langle X\rangle$ such that
$\delta_0(x_j)=f_j+T({\mathfrak W})\in F_{\infty}({\mathfrak W})$, $j=1,2,\ldots$.
Since every mapping $\{x_1,x_2,\ldots\}\to K\langle x_1,x_2,\ldots\rangle$
can be extended to a derivation of $K\langle x_1,x_2,\ldots\rangle$,
it is sufficient to show that the derivation $\Delta$ of
$K\langle x_1,x_2,\ldots\rangle$ defined by $\Delta(x_j)=f_j$, $j=1,2,\ldots$,
has the property $\Delta(T({\mathfrak W}))\subset T({\mathfrak W})$.
Since the field $K$ is of characteristic 0, if
$u(x_1,\ldots,x_m)$ belongs to $T({\mathfrak W})$,
then the multihomogeneous components
of $u$ also are in $T({\mathfrak W})$ and we may assume that
$u(x_1,\ldots,x_m)\in T({\mathfrak W})$ is multihomogeneous.
The partial linearization $u_j(x_1,\ldots,x_m,x_{m+1})$
in $x_j$ of $u(x_1,\ldots,x_m)$, i.e. the linear component
in $x_{m+1}$ of $u(x_1,\ldots,x_{j-1},x_j+x_{m+1},x_{j+1},\ldots,x_m)$
also belongs to $T({\mathfrak W})$. It is easy to see that $\Delta$ acts on
$u(x_1,\ldots,x_m)$ by
\[
\Delta(u(x_1,\ldots,x_m))=\sum_{j=1}^mu_j(x_1,\ldots,x_m,\Delta(x_j)).
\]
Since $u_j(x_1,\ldots,x_m,\Delta(x_j))\in T({\mathfrak W})$ we obtain that
$\Delta(u)\in T({\mathfrak W})$ and this means that $\Delta$ induces a derivation
$\delta$ on
$F_{\infty}({\mathfrak W})=K\langle x_1,x_2,\ldots\rangle/T({\mathfrak W})$
with the additional property $\delta(x_j)=f_j$, and $\delta$ extends $\delta_0$.
This implies also the case of $F_m({\mathfrak W})$:
If $f_1,\ldots,f_m\in F_m({\mathfrak W})$, then we extend the mapping to
a derivation of $F_{\infty}({\mathfrak W})$ (e.g. by $\delta_0(x_j)=0$ for $j>m$).
Then the restriction to $F_m({\mathfrak W})$ of the derivation of
$F_{\infty}({\mathfrak W})$ is a derivation of $F_m({\mathfrak W})$.
\end{proof}

\section{Weitzenb\"ock Derivations of Polynomial Algebras}

Since we consider nonzero Weitzenb\"ock derivations only, without loss
of generality we may assume that the derivation $\delta$ is
in its Jordan normal form,
$\delta(x_1)=0$, $\delta(x_2)=x_1$ and the set of variables
$X=\{x_1,\ldots,x_m\}$ is a Jordan basis of $V_m=\sum_{j=1}^mKx_j$.
If the rank of $\delta$ is equal to $m-1$ (i.e. $\delta(x_j)=x_{j-1}$,
$j=2,\ldots,m$), then $\delta$ is called the basic Weitzenb\"ock derivation
of $K[X]$. The following proposition, see \cite{No}, gives the description of
the algebras of constants of any Weitzenb\"ock derivation. (It is a very special case
of the more general situation of an arbitrary locally nilpotent derivation.)
For our purposes we work in the localization of the polynomial algebra
$K[X][x_1^{-1}]=K[x_1,x_2,\ldots,x_m][x_1^{-1}]$ consisting
of all polynomials in $x_1,\ldots,x_m$ allowing negative degrees of $x_1$.
Since $x_1$ is a constant (i.e. $\delta(x_1)=0$), we may extend $\delta$ to
a derivation of $K[X][x_1^{-1}]$.

\begin{proposition} \label{generators Weitzenboeck}
Let $\delta^{p_j+1}(x_j)=0$, $j=1,\ldots,m$, and let
\[
z_j=\sum_{k=0}^{p_j}\frac{\delta^k(x_j)}{k!}(-x_2)^kx_1^{p_j-k},\quad
j=3,4,\ldots,m,
\].
\par
{\rm (i)} $(K[X][x_1^{-1}])^{\delta}=K[x_1,z_3,z_4,\ldots,z_m][x_1^{-1}]$;
\par
{\rm (ii)} $K[X]^{\delta}=K[X]\cap(K[X][x_1^{-1}])^{\delta}$.
\end{proposition}

\begin{example}
If $\delta$ is a basic Weitzeb\"ock derivation, then
\[
z_3=x_3x_1^2-\frac{x_2^2x_1}{2}=\frac{x_1}{2}(2x_3x_1-x_2^2),\quad
z_4=x_1\left(x_4x_1^2-x_3x_2x_1+\frac{x_2^3}{3}\right),\ldots,
\]
\[
z_j=(-1)^j\frac{j}{(j+1)!}x_1\left(x_2^{j+1}
+\frac{(j+1)!}{j}\sum_{k=0}^{j-1}(-1)^{j-k}\frac{1}{(j+1)!}x_1^{j-k}x_2^kx_{j+1-k}\right).
\]
\end{example}

\begin{corollary} \label{transcendence degree}
For any Weitzenb\"ock derivation $\delta$,
the transcendence degree (i.e. the maximal number of algebraically independent
elements) of $K[x_1,\ldots,x_m]^{\delta}$ is equal to $m-1$.
\end{corollary}

The explicit form of the generators of $K[x_1,\ldots,x_m]^{\delta}$
is known for small $m$ only. Tan \cite{Ta} presented an algorithm for
computing the generators of the algebra of constants of a basic derivation.
It was generalized by van den Essen \cite{E1} for any locally nilpotent derivation
assuming that the finite generation of the algebra of constants is known.
The algorithm involves Gr\"obner bases techniques.

\begin{examples} \label{examples of concrete generators}
We have selected few examples of the generating sets of the algebra of constants,
all of them taken from \cite{No}.
For $\delta$ being a basic Weitzenb\"ock derivation
(with $\delta(x_1)=0$ and $\delta(x_j)=x_{j-1}$, $j=2,\ldots,m$):
\[
K[x_1,x_2]^{\delta}=K[x_1],\quad
K[x_1,x_2,x_3]^{\delta}=K[x_1,x_2^2-2x_1x_3],
\]
\[
K[x_1,x_2,x_3,x_4]^{\delta}=K[x_1,x_2^2-2x_1x_3, x_2^3-3x_1x_2x_3+3x_1^2x_4,
\]
\[
x_2^2x_3^2-2x_2^3x_4+6x_1x_2x_3x_4-\frac{8}{3}x_1x_3^3-3x_1^2x_4^2],
\]
(see \cite{No}, Example 6.8.2),
\[
K[x_1,x_2,x_3,x_4,x_5]^{\delta}=K[x_1,x_2^2-2x_1x_3,
2x_2x_4-x_3^2-2x_1x_5,
\]
\[
x_2^3-3x_1x_2x_3+3x_1^2x_4,
6x_2^2x_5-6x_2x_3x_4+2x_3^3-12x_1x_3x_5+9x_1x_4^2],
\]
(see \cite{No}, Example 6.8.4).
\par
For $\delta$ nonbasic, $\delta(x_2)=x_1$, $\delta(x_4)=x_3$,
$\delta(x_1)=\delta(x_3)=0$ (see \cite{No}, Proposition 6.9.5):
\[
K[x_1,x_2,x_3,x_4]^{\delta}=K[x_1,x_3, x_1x_4-x_2x_3],
\]
for $\delta$ defined by $\delta(x_3)=x_2$, $\delta(x_2)=x_1$,
$\delta(x_5)=x_4$,
$\delta(x_1)=\delta(x_4)=0$ (see \cite{No}, Example 6.8.5):
\[
K[x_1,x_2,x_3,x_4,x_5]^{\delta}=K[x_1,x_4, x_1x_5-x_2x_4,x_2^2-2x_1x_3,
2x_3x_4^2-2x_2x_4x_5+x_1x_5^2].
\]
\end{examples}

\begin{remark}
Springer \cite{Sr} found a formula for the Hilbert series of the algebra of
invariants of $SL_2(K)$ acting on the forms of degree $d$. This is equivalent
to the description of the Hilbert series of the algebra of constants of the
basic Weitzenb\"ock derivation of $K[x_1,\ldots,x_{d+1}]$. Almkvist
\cite{A1, A2} related these invariants with invariants of the modular action
of a cyclic group of order $p$.
\end{remark}

\section{Lifting and Description of the Constants}

We need the following easy lemma.

\begin{lemma} \label{lifting in general case}
Let $G\subset H$ be groups and let the $H$-module $M$ be completely
reducible. If $N\subset M$ is an $H$-submodule and $\bar m\in M/N$
is a $G$-invariant, then $\bar m$ can be lifted to a $G$-invariant $m\in M$.
\end{lemma}

\begin{proof}
Let $P$ be an $H$-complement of $N$ in $M$, i.e. $M=N\oplus P$. Since
$M/N\cong P$, there exists an element $m\in P$ which maps on $\bar m$
under the natural homomorphism $M\to M/N$.
Since $\bar m$ is $G$-invariant, we obtain that $\overline{G(m)}=G(\bar m)=\bar m$.
Taking into account that $m_1,m_2\in P$, $m_1\not= m_2$,
implies that $\bar m_1\not= \bar m_2$ in $M/N$, and $G(P)=P$,
we deduce that $G(m)=m$ in $M$, i.e. $m$ is $G$-invariant.
\end{proof}

\begin{proposition} \label{lifting of invariants of linear groups}
Let $K\langle X\rangle=K\langle x_1,\ldots,x_m\rangle$ be the free
associative algebra with the canonical $GL_m$-action, and let
$I\subset J$ be $GL_m$-invariant two-sided ideals of $K\langle X\rangle$.
Then for every subgroup $G$ of $GL_m$, the $G$-invariants of
$K\langle X\rangle/J$ can be lifted to $G$-invariants of $K\langle X\rangle/I$.
\end{proposition}

\begin{proof} The statement follows immediately from
Lemma \ref{lifting in general case} because, as a $GL_m$-module,
$K\langle X\rangle$ is completely reducible.
\end{proof}

\begin{corollary} \label{lifting of constants}
If $I\subset J$ are $GL_m$-invariant two-sided ideals of
$K\langle X\rangle$ and $\delta$ is
a Weitzenb\"ock derivation on $K\langle X\rangle$, then the algebra
of constants $(K\langle X\rangle/J)^{\delta}$
can be lifted to the algebra of constants $(K\langle X\rangle/I)^{\delta}$.
\end{corollary}

\begin{proof} The corollary is a straightforward consequence of
Proposition \ref{lifting of invariants of linear groups} because
the algebras of constants $(K\langle X\rangle/J)^{\delta}$
and $(K\langle X\rangle/I)^{\delta}$ coincide, respectively,
with the algebras of $g$-invariants
$(K\langle X\rangle/J)^g$ and $(K\langle X\rangle/I)^g$, where
$g=\exp\delta$ is the linear transformation corresponding to $\delta$.
\end{proof}

\begin{corollary} \label{test for finite generation}
Let $I\subset J$ be $GL_m$-invariant two-sided ideals of
$K\langle X\rangle$ and let $\delta$ be
a Weitzenb\"ock derivation on $K\langle X\rangle$. If the algebra
of constants $(K\langle X\rangle/J)^{\delta}$ is not finitely generated, then
$(K\langle X\rangle/I)^{\delta}$ is also not finitely generated.
\end{corollary}

\begin{remark}
Corollary \ref{test for finite generation} holds also for Lie algebras and other
free algebras including free (special or not) Jordan algebras and the absolutely
free algebra $K\{x_1,\ldots,x_m\}$.
\end{remark}

Now we shall describe the algebras of constants in the case of two variables,
assuming that $K\langle x_1,x_2\rangle=K\langle x,y\rangle$ and
$\delta(x)=0$, $\delta(y)=x$.
\par
Recall that any irreducible polynomial $GL_2$-module $W(\lambda_1,\lambda_2)$
has a unique (up to a multiplicative constant) element $w(x,y)$ which is
bihomogeneous of degree $(\lambda_1,\lambda_2)$ and is called the highest weight
vector of $W(\lambda_1,\lambda_2)$. For any $GL_2$-invariant homomorphic image
$K\langle x,y\rangle/I$ of $K\langle x,y\rangle$ the algebra of constants
$(K\langle x,y\rangle/I)^{\delta}$ coincides with the algebra of
$g$-invariants $(K\langle x,y\rangle/I)^g$ where $g=\exp\delta$. Since
$g(x)=x$, $g(y)=x+y$ and $\text{\rm char}K=0$, the algebra of $g$-invariants
coincides with the algebra of invariants of the unitriangular group $UT_2(K)$.
Hence, as in Almkvist, Dicks and Formanek \cite{ADF}, we may use Theorem 3.3 (i)
of De Concini, Eisenbud and Procesi \cite{DEP} and obtain:

\begin{theorem} \label{constants in two variables}
For any $GL_2$-invariant ideal $I$ of $K\langle x,y\rangle$ the algebra
of constants $(K\langle x,y\rangle/I)^{\delta}$ is spanned by the highest weight
vectors of the $GL_2$-irreducible components of $K\langle x,y\rangle/I$.
\end{theorem}

\begin{remarks} \label{three remarks}
1. A direct proof of Theorem \ref{constants in two variables} can be obtained
using the criterion of Koshlukov \cite{Ko1} which states:
A multihomogeneous of degree $\lambda=(\lambda_1,\ldots,\lambda_m)$ polynomial
$w(x_1,\ldots,x_m)\in K\langle x_1,\ldots,x_m\rangle$
is a highest weight vector of an irreducible $GL_m$-submodule
$W(\lambda)$ of $K\langle x_1,\ldots,x_m\rangle$ if and only if
for all partial linearizations
$w_j(x_1,\ldots,x_m,x_{m+1})$ of $w(x_1,\ldots,x_m)$ one has
$w_j(x_1,\ldots,x_m,x_i)=0$ for all $i<j$.
\par
2. By Almkvist, Dicks and Formanek \cite{ADF} the algebra
$(K\langle x_1,\ldots,x_m\rangle)^{UT_m(K)}$ of all
$UT_m(K)$-invariants coincides with the vector space spanned by all
highest weight vectors $w(x_1,\ldots,x_m)\in W(\lambda)\subset
K\langle x_1,\ldots,x_m\rangle$,
when $\lambda=(\lambda_1,\ldots,\lambda_m)$
runs on the set of all partitions in not more than $m$ parts.
\par
3. Following Almkvist, Dicks and Formanek \cite{ADF},
for any unipotent transformation $g$ of $K\langle x_1,\ldots,x_m\rangle$
(and hence for any Weitzenb\"ock derivation $\delta$)
one can define a $GL_2$-action on $K\langle x_1,\ldots,x_m\rangle$
and on the factor algebras $K\langle x_1,\ldots,x_m\rangle/I$
modulo $GL_m$-invariant ideals,
such that $(K\langle x_1,\ldots,x_m\rangle)^g$ and
$(K\langle x_1,\ldots,x_m\rangle/I)^g$ are spanned by the highest weight vectors
with respect to the $GL_2$-action.
\end{remarks}

The necessary background on symmetric functions which we need  can be found
e.g. in the book by Macdonald \cite{M}.
Any symmetric function in $m$ variables $f(t_1,\ldots,t_m)$
which can be expressed as a formal power series has the presentation
\[
f(t_1,\ldots,t_m)=\sum_{\lambda}m(\lambda)S_{\lambda}(t_1,\ldots,t_m),
\]
where $S_{\lambda}(t_1,\ldots,t_m)$ is the Schur function corresponding
to the partition $\lambda=(\lambda_1,\ldots,\lambda_m)$ and $m(\lambda)$
is the multiplicity of $S_{\lambda}(t_1,\ldots,t_m)$ in $f(t_1,\ldots,t_m)$.
This presentation agrees with the theory of polynomial representations of $GL_m$
because the Schur functions play the role of characters of the irreducible
polynomial $GL_m$-representations. In our case this relation gives the following:
If $K\langle X\rangle/I$ for some $GL_m$-invariant ideal $I$, then the Hilbert series
of $K\langle X\rangle/I$ has the presentation
\[
H(K\langle X\rangle/I,t_1,\ldots,t_m)
=\sum_{\lambda}m(\lambda)S_{\lambda}(t_1,\ldots,t_m),
\]
if and only if $K\langle X\rangle/I$ is decomposed as a $GL_m$-module as
\[
K\langle X\rangle/I\cong\sum_{\lambda}m(\lambda)W(\lambda).
\]
In the case of two variables the Schur functions have the following simple expression
\[
S_{(\lambda_1,\lambda_2)}(t_1,t_2)=(t_1t_2)^{\lambda_2}
\frac{t_1^{\lambda_1-\lambda_2+1}-t_2^{\lambda_1-\lambda_2+1}}{t_1-t_2}.
\]
Drensky and Genov \cite{DGn} defined the multiplicity series of
\[f(t_1,t_2)=\sum_{\lambda}m(\lambda)S_{\lambda}(t_1,t_2)
\]
as the formal power series
\[
M(f)(t,u)=\sum_{\lambda}m(\lambda)t^{\lambda_1}u^{\lambda_2},
\]
or, if one introduces a new variable $v=tu$, as
\[
M'(f)(t,v)=\sum_{\lambda}m(\lambda)t^{\lambda_1-\lambda_2}v^{\lambda_2}.
\]
The relation between the symmetric function and its multiplicity series is
\[
f(t_1,t_2)=\frac{t_1M'(f)(t_1,t_1t_2)-t_2M'(f)(t_2,t_1t_2)}{t_1-t_2}.
\]
Theorem \ref{constants in two variables} gives that the Hilbert series
of the algebra of constants $(K\langle x,y\rangle/I)^{\delta}$
(with respect to the bigrading) is equal to
the multiplicity series of the Hilbert series of $K\langle x,y\rangle/I$:

\begin{corollary} \label{Hilbert series of constants}
For any $GL_2$-invariant ideal $I$ of $K\langle x,y\rangle$
and for the basic Weitzen\-b\"ock derivation $\delta$
\[
H((K\langle x,y\rangle/I)^{\delta},t,u)
=M(H(K\langle x,y\rangle/I)(t,u).
\]
\end{corollary}

If we consider the usual grading, Corollary
\ref{Hilbert series of constants}
has the form
\[
H((K\langle x,y\rangle/I)^{\delta},t)
=M(H(K\langle x,y\rangle/I)(t,t)
=M'(H(K\langle x,y\rangle/I)(t,t^2).
\]
We shall apply Corollary
\ref{Hilbert series of constants}
in the next section in the concrete description
of the generators of the constants in $K\langle x,y\rangle$ and, more generally,
in any relatively free associative algebra.

\section{Examples and Concrete Generators of Algebras of Constants}

We start this section with several examples when we determine completely
the algebras of constants and their generators. We shall consider algebras of rank
2 and 3 only and shall denote the free generators by $x,y$ and $x,y,z$,
respectively.
We shall handle the case of basic Weitzenb\"ock derivations $\delta$ only,
assuming that $\delta(x)=0$, $\delta(y)=x$
(and $\delta(z)=y$ if the rank of the algebra is equal to 3).

\begin{example} \label{Grassmann}
Let ${\mathfrak L}_2$ be the variety of associative algebras defined
by the identity $[[x,y],z]=0$. By the theorem of Krakowski and Regev \cite{KR}
${\mathfrak L}_2$ coincides with the variety generated by the infinite dimensional
Grassmann algebra. The $S_n$-cocharacter sequence of ${\mathfrak L}_2$ is
equal to
\[
\chi_n({\mathfrak L}_2)=\sum_{k=1}^n\chi_{(k,1^{n-k})},\quad n\geq 1,
\]
see \cite{KR}. In virtue of the correspondence between
cocharacters and Hilbert series, see \cite{Bl} and \cite{D1}
(or the book \cite{D6}) the Hilbert series of the relatively free algebra
$F_m({\mathfrak L}_2)$ is equal to
\[
H(F_m({\mathfrak L}_2),t_1,\ldots,t_m)
=1+\sum_{k\geq 1}\sum_{l=0}^{m-1}S_{(k,1^l)}(t_1,\ldots,t_m).
\]
It is well known that $F_m({\mathfrak L}_2)$ has a basis
\[
x_1^{a_1}\cdots x_m^{a_m}[x_{i_1},x_{i_2}]\cdots[x_{i_{2p-1}},x_{i_{2p}}],
\quad 1\leq i_1<i_2<\cdots<i_{2p-1}<i_{2p}\leq m,
\]
see for example Bokut and Makar-Limanov \cite{BM} or the book \cite{D6}.
The commutators $[x_i,x_j]$ are in the centre of $F_m({\mathfrak L}_2)$
and satisfy the relations
\[
[x_{\sigma(1)},x_{\sigma(2)}]\cdots[x_{\sigma(2p-1)},x_{\sigma(2p)}]
=(\text{\rm sign}\sigma)[x_1,x_2]\cdots[x_{2p-1},x_{2p}],\quad \sigma\in S_{2p}.
\]
\par
Let $m=2$. Then $F_2({\mathfrak L}_2)$ has a basis
\[
\left\{ x^ay^b,x^ay^b[x,y]\mid a,b\geq 0\right\}.
\]
Its Hilbert series and the related multiplicity series are, respectively,
\[
H(F_2({\mathfrak L}_2),t_1,t_2)
=\frac{1+t_1t_2}{(1-t_1)(1-t_2)}
=\sum_{n\geq 0}S_{(n)}(t_1,t_2)
+\sum_{n\geq 2}S_{(n-1,1)}(t_1,t_2),
\]
\[
M(H(F_2({\mathfrak L}_2))(t,u)=\sum_{n\geq 0}t^n+\sum_{n\geq 2}t^{n-1}u
=\frac{1+tu}{1-t}.
\]
By Corollary \ref{Hilbert series of constants},
\[
H((F_2({\mathfrak L}_2))^{\delta},t,u)=\frac{1+tu}{1-t}.
\]
Since the vector subspace of $F_2({\mathfrak L}_2)$ spanned by
$x^n$, $n\geq 0$, and $x^{n-2}[x,y]$, $n\geq 2$, consists of $\delta$-constants
and has the same Hilbert series as $(F_2({\mathfrak L}_2))^{\delta}$,
we obtain that it coincides with the algebra of constants. This immediately
implies that the algebra $(F_2({\mathfrak L}_2))^{\delta}$ is generated by $x$
and $[x,y]$.
\par
Let $m=3$. Then $F_3({\mathfrak L}_2)$ has a basis
\[
\left\{x^ay^bz^c,x^ay^bz^c[x,y],x^ay^bz^c[x,z],x^ay^bz^c[y,z]
\mid a,b,c\geq 0 \right\}
\]
and the commutator ideal $C$ of $F_3({\mathfrak L}_2)$ is a free $K[x,y,z]$-module
with free generators $[x,y],[x,z],[y,z]$.
By Examples \ref{examples of concrete generators},
$K[x,y,z]^{\delta}=K[x,y^2-2xz]$. We may choose $y^2-xz-zx$ as
a lifting in $(F_3({\mathfrak L}_2))^{\delta}$
of $y^2-2xz\in (K[x,y,z])^{\delta}$. Hence $(F_3({\mathfrak L}_2))^{\delta}$
is generated by $x,y^2-xz-zx$ and some elements in the commutator ideal $C$.
Every element of $K[x,y,z]$ can be written in a unique way as
\[
f_0(x,y^2-2xz)+\sum_{n\geq 1}f_n(x,y^2-2xz)z^n+\sum_{n\geq 1}g_n(x,y^2-2xz)yz^{n-1}.
\]
Hence the elements in $C$ have the form
\[
f=\alpha(x,y,z)[x,y]+\beta(x,y,z)[x,z]+\gamma(x,y,z)[y,z],\quad
\alpha,\beta,\gamma\in K[x,y,z].
\]
If $f$ is a $\delta$-constant, then
\[
0=\delta(f)=(\delta(\alpha)+\beta)[x,y]+(\delta(\beta)+\gamma)[x,z]+\delta(\gamma)[y,z].
\]
In this way, $f\in (F_3({\mathfrak L}_2))^{\delta}$
if and only if
\[
\delta(\gamma)=0,\quad \delta(\beta)=-\gamma,\quad
\delta(\alpha)=-\beta.
\]
We present $\beta(x,y,z)$ in the form
\[
\beta=f_0+\sum_{n\geq 1}(f_nz^n+g_nz^{n-1}y),
\quad f_0,f_n,g_n\in (K[x,y,z])^{\delta},
\]
and calculate, bearing in mind that $y^2=(y^2-2xz)+2xz$,
\[
-\gamma=\delta(\beta)
=\sum_{n\geq 1}\left(nf_nz^{n-1}y+(n-1)g_nz^{n-2}y^2+xg_nz^{n-1}\right)
\]
\[
=\sum_{n\geq 1}\left((n-1)g_n(y^2-2xz)z^{n-2}
+(2n-1)xg_nz^{n-1}+nf_nz^{n-1}y\right).
\]
This easily implies that $f_n=0$, $n\geq 1$, $g_n=0$, $n\geq 2$, and
$\beta=f_0+g_1y$, $f_0,g_1\in K[x,y^2-2xz]=(K[x,y,z])^{\delta}$. Hence
$\gamma=-g_1x$. Continuing in this way, we obtain the final form of $\alpha,\beta,\gamma$:
\[
\alpha=\alpha_0z+\alpha_1y+\alpha_2,\quad \beta=-\alpha_0y-\alpha_1x,\quad \gamma=\alpha_0x.
\]
Hence the part of the algebra of constants of $F_3({\mathfrak L}_2)$ which belongs
to the commutator ideal $C$ is spanned as a $(K[x,y,z])^{\delta}$-module by
\[
[x,y], \quad y[x,y]-x[x,z],\quad z[x,y]-y[x,z]+x[y,z],
\]
and $(F_3({\mathfrak L}_2))^{\delta}$ is generated by
\[
x,\quad y^2-xz-zx,\quad
[x,y], \quad y[x,y]-x[x,z],\quad z[x,y]-y[x,z]+x[y,z].
\]
\end{example}

\begin{example} \label{free metabelian associative algebra}
Let us consider the variety $\mathfrak M$ of all ``metabelian'' associative
algebras defined by the identity $[x_1,x_2][x_3,x_4]=0$. It is well known
that $F_2({\mathfrak M})$ has a basis
\[
\{x^ay^b,\quad x^ay^b[x,y]x^cy^d\mid a,b,c,d\geq 0\}.
\]
We shall write the element $x^ay^b[x,y]x^cy^d$ as
$[x,y]x_1^ay_1^bx_2^cy_2^d$. In this way,
the commutator ideal $C$ of $F_2({\mathfrak M})$ is a free
cyclic $K[x,y]$-bimodule (or a free
cyclic $K[x_1,y_1,x_2,y_2]$-module) with the $K[x,y]$-action defined by
\[
x[x,y]=[x,y]x_1,\quad y[x,y]=[x,y]y_1,\quad
[x,y]x=[x,y]x_2,\quad [x,y]y=[x,y]y_2.
\]
The Hilbert series of $F_2({\mathfrak M})$ is
\[
H(F_2({\mathfrak M}),t_1,t_2)
=\frac{1}{(1-t_1)(1-t_2)}+\frac{t_1t_2}{(1-t_1)^2(1-t_2)^2}.
\]
One can calculate directly the $S_n$-cocharacter of $\mathfrak M$
using the Young rule as in \cite{D6} or to apply techniques of \cite{DGn}
to see that the multiplicty series of $F_2({\mathfrak M})$ is
\[
M'(F_2({\mathfrak M}))(t,v)=\frac{1}{1-t}+\frac{v}{(1-t)^2(1-v)}.
\]
By Corollary \ref{Hilbert series of constants} this is also the Hilbert series of
the algebra of constants $(F_2({\mathfrak M}))^{\delta}$.
We consider the linearly independent highest weight vectors
\[
x^n,\ n\geq 0,\quad
[x,y]x_1^px_2^q(x_1y_2-y_1x_2)^r,\ p,q,r\geq 0.
\]
They span a graded vector subspace of $(F_2({\mathfrak M}))^{\delta}$
and its Hilbert series coincides with the HIlbert series of
$(F_2({\mathfrak M}))^{\delta}$. Hence the above highest weight vectors span
$(F_2({\mathfrak M}))^{\delta}$. Since the square of
the commutator ideal $C$ is equal to 0,
the element $x$ together with all
$[x,y]x_1^px_2^q(x_1y_2-y_1x_2)^r$, $p,q,r\geq 0$, is a minimal generating set
of $(F_2({\mathfrak M}))^{\delta}$
and the algebra of constants is not finitely generated.
\end{example}

Now we start with the description of the constants of the free algebra
$K\langle x,y\rangle$ which will gives also the description of
the constants in any two-generated associative algebra.

\begin{proposition} \label{Hilbert series of constants of free associative algebra}
The Hilbert series of the algebra of constants $(K\langle x,y\rangle)^{\delta}$ are
\[
H((K\langle x,y\rangle)^{\delta},t,u)=\sum_{(\lambda_1,\lambda_2)}
\left(\binom{\lambda_1+\lambda_2}{\lambda_2}
-\binom{\lambda_1+\lambda_2}{\lambda_2-1}\right)t^{\lambda_1}u^{\lambda_2}
\]
\[
=\frac{1-\sqrt{1-4v}}{2v}\cdot
\frac{1}{1-\frac{1-\sqrt{1-4v}}{2v}t},
\]
where $v=tu$ and, in one variable,
\[
H((K\langle x,y\rangle)^{\delta},t)=\sum_{p\geq 0}
\left(\binom{2p}{p}t^{2p}+\binom{2p+1}{p}t^{2p+1}\right).
\]
\end{proposition}

\begin{proof}
By Corollary \ref{Hilbert series of constants} the Hilbert series of
the algebra of constants $(K\langle x,y\rangle)^{\delta}$
is equal to the multiplicity series of
the Hilbert series of $K\langle x,y\rangle$.
By representation theory of general linear groups, the multiplicity $m_{\lambda}$
of the irreducible $GL_m$-module $W(\lambda)$ in $K\langle x_1,\ldots,x_m\rangle$
for the partition $\lambda$ of $n$
is equal to the degree $d_{\lambda}$
of the irreducible $S_n$-character $\chi_{\lambda}$.
By the hook formula, for $\lambda=(\lambda_1,\lambda_2)$
\[
d_{\lambda}=\frac{(\lambda_1+\lambda_2)(\lambda_1-\lambda_2+1)}{(\lambda_1+1)!\lambda_2!}
=\binom{\lambda_1+\lambda_2}{\lambda_2}-\binom{\lambda_1+\lambda_2}{\lambda_2-1}.
\]
This gives the expression for
$H((K\langle x,y\rangle)^{\delta},t,u)$. If we set there $u=t$ we obtain
that the coefficient of $t^{2p}$ is equal to
\[
\sum_{i=0}^p\left(\binom{2p}{i}-\binom{2p}{i-1}\right)=\binom{2p}{p},
\]
and similarly for the coefficient of $t^{2p+1}$.
In order to obtain the formula in terms of $t$ and $v$ we can either use the
known formulas for the summation of formal power series with binomial coefficients
or proceed in the following way using ideas from \cite{DGn}.
The Hilbert series of $K\langle x,y\rangle$ is equal to
\[
f(t_1,t_2)=H(K\langle x,y\rangle,t_1,t_2)=\frac{1}{1-(t_1+t_2)}.
\]
It is sufficient to show that the multiplicity series of $f(t_1,t_2)$ is
\[
M'(f)(t,v)=\frac{1-\sqrt{1-4v}}{2v}\cdot
\frac{1}{1-\frac{1-\sqrt{1-4v}}{2v}t}.
\]
Since the multiplicity series of any symmetric function
$f(t_1,t_2)\in K[[t_1,t_2]]$ is a uniquely determined formal power series
in $K[[t,v]]$, it is sufficient to show that the expansion of
\[
\frac{1-\sqrt{1-4v}}{2v}\cdot
\frac{1}{1-\frac{1-\sqrt{1-4v}}{2v}t}
\]
is in $K[[t,v]]$ (which is obvious because $1-\sqrt{1-4v}=\sum_{n\geq 1}a_nv^n$
for some $a_n\in K$ and $\frac{1-\sqrt{1-4v}}{2v}\in K[[v]]$) and
to use the formula
\[
\frac{t_1M'(f)(t_1,t_1t_2)-t_2M'(f)(t_2,t_1t_2)}{t_1-t_2}=f(t_1,t_2).
\]
Direct verification shows that for
\[
g(t,v)=\frac{1-\sqrt{1-4v}}{2v}\cdot
\frac{1}{1-\frac{1-\sqrt{1-4v}}{2v}t}
\]
\[
\frac{t_1g(t_1,t_1t_2)-t_2g(t_2,t_1t_2)}{t_1-t_2}=\frac{1}{1-(t_1+t_2)}
\]
which gives that $g(t,v)=M'(f)(t,v)$.
\end{proof}

By the theorem of Lane \cite{L} and Kharchenko \cite{K1},
the algebra of constants $(K\langle X\rangle)^{\delta}$ is a graded free algebra
and hence has a homogeneous system of free generators. The following theorem
describes the generating function of the set of free generators.

\begin{theorem}\label{generating function of the algebra of constants}
The generating function of any bihomogeneous system of free generators of
$(K\langle X\rangle)^{\delta}$ with respect to the variables $t$ and $v=tu$ is
\[
a(t,v)=t+\frac{1-\sqrt{1-4v}}{2}.
\]
\end{theorem}

\begin{proof}
If $a(t,v)$ is the generating function of the set of free generators
of $(K\langle X\rangle)^{\delta}$, then the Hilbert series of
$(K\langle X\rangle)^{\delta}$ is
\[
H((K\langle X\rangle)^{\delta},t,v)=\frac{1}{1-a(t,v)}.
\]
Applying Proposition \ref{Hilbert series of constants of free associative algebra}
we obtain that
\[
\frac{1}{1-a(t,v)}
=\frac{1-\sqrt{1-4v}}{2v}\cdot
\frac{1}{1-\frac{1-\sqrt{1-4v}}{2v}t}
\]
and the expression of $a(t,v)$ is a result of easy calculations.
\end{proof}

\begin{corollary} \label{generatation by x and sl-invariants}
The algebra of constants $(K\langle x,y\rangle)^{\delta}$, where
$\delta(x)=0$, $\delta(y)=x$, is generated by $x$ and by $SL_2(K)$-invariants.
\end{corollary}

\begin{proof}
An element $f(x,y)\in K\langle x,y\rangle$ is a $SL_2$-invariant
if and only if it is a linear combination of
highest weight vectors of a $GL_2$-submodules
$W(\lambda_1,\lambda_1)$. By Theorem \ref{constants in two variables},
the $\delta$-constants are linear combinations of highest weight vectors
$w_{(\lambda_1,\lambda_2)}$, and  $w_{(\lambda_1,\lambda_2)}$
is bihomogeneous of degree $(\lambda_1,\lambda_2)$. Hence we obtain that
the set of $SL_2$-invariants coincides with the linear combinations of
bihomogeneous elements of degree $(p,p)$.
The only nonzero coefficients of the Hilbert series
$H((K\langle x,y\rangle)^{SL_2},t,u)$ are for $v^n=(tu)^n$ and
$H((K\langle x,y\rangle)^{SL_2},t,u)$ is obtained from
$H((K\langle x,y\rangle)^{\delta},t,v)$ by replacing $t$ with 0.
Hence Theorem \ref{generating function of the algebra of constants}
gives that the set of homogeneous generators
of the algebra of $\delta$-constants is spanned by $x$ and $SL_2$-invariants.
\end{proof}

Corollary \ref{lifting of constants} gives immediately:

\begin{corollary} \label{generatation by sl-invariants for algebras of rank two}
For any $GL_2$-invariant ideal $I$ of $K\langle x,y\rangle$
the algebra of constants $(K\langle x,y\rangle/I)^{\delta}$, where
$\delta(x)=0$, $\delta(y)=x$, is generated by $x$ and by $SL_2(K)$-invariants.
\end{corollary}

\begin{remark}\label{Catalan}
By Almkvist, Dicks and Formanek \cite{ADF} Example 5.10,
the Hilbert series of the algebra
of $SL_2$-invariants of $K\langle x,y\rangle$ is
\[
H((K\langle x,y\rangle)^{SL_2},v)=\frac{1-\sqrt{1-4v}}{2v}
=\sum_{n\geq 0}\frac{1}{n+1}\binom{2n}{n}v^n,
\]
and the coefficient of $v^n$ is the $(n+1)$-st Catalan number
$c_{n+1}$. (By definition $c_n$ is
the number of possibilities to distribute parentheses in the sum
$1+1+\cdots+1$ of $n$ units, see e.g. \cite{Ha}.)
This agrees with Proposition
\ref{Hilbert series of constants of free associative algebra} because
$H((K\langle x,y\rangle)^{SL_2},v)$ is obtained from
\[
H((K\langle x,y\rangle)^{\delta},t,v)
=\frac{1-\sqrt{1-4v}}{2v}\cdot
\frac{1}{1-\frac{1-\sqrt{1-4v}}{2v}t}
\]
by replacing $t$ with 0.
\par
Theorem \ref{generating function of the algebra of constants}
gives that the generating function of a homogeneous system of free generators
of $(K\langle X\rangle)^{SL_2}$ is
\[
b(v)=\frac{1-\sqrt{1-4v}}{2}=vH((K\langle x,y\rangle)^{SL_2},v).
\]
Since $v=tu$ is of second degree, the number of generators
of $(K\langle x,y\rangle)^{SL_2}$ of degree $2n$ is equal to the $n$-th
Catalan number.
\end{remark}

Below we give an inductive procedure to construct
an infinite set of free generators
of the algebra $(K\langle x,y\rangle)^{SL_2}$.

\begin{algorithm}
The following infinite procedure gives a complete set
$\{w_1,w_2,\ldots\}$ of free
generators of the algebra $(K\langle x,y\rangle)^{SL_2}$.
We set $w_1=[x,y]$. If we have already constructed all
free generators $w_1,w_2,\ldots,w_k$ of degree $\leq 2n$, then
we form all $c_{n+1}$ products
$w_{i_1}\cdots w_{i_s}$ of degree $2n$, which we number as
$\omega_j$, $j=1,\ldots,c_{n+1}$,
and add to the system of generators the $c_{n+1}$
elements
\[
w_{k+j}=x\omega_jy-y\omega_jx
=xw_{i_1}\cdots w_{i_s}y-yw_{i_1}\cdots w_{i_s}x,\quad
j=1,\ldots,c_{n+1}.
\]
\end{algorithm}

The first several elements of the generating set are:
\[
w_1=[x,y],\quad w_2=x[x,y]y-y[x,y]x,
\]
\[
w_3=xw_1^2y-yw_1^2x=x[x,y]^2y-y[x,y]^2x,
\]
\[
w_4=xw_2y-yw_2x=x(x[x,y]y-y[x,y]x)y-y(x[x,y]y-y[x,y]x)x.
\]

\begin{proof}
By Remark \ref{Catalan} and by inductive arguments, we may assume that
the number of products $\omega_j=w_{i_1}\cdots w_{i_s}$ of degree $2n$ is equal
to the Catalan number $c_{n+1}$. Hence the number of words
$x\omega_jy-y\omega_jx$, all of degree $2(n+1)$ is also equal to $c_{n+1}$
which agrees with the number of free generators of degree $2(n+1)$.
Clearly, if $\omega_j$ is an $SL_2$-invariant, the element
$x\omega_jy-y\omega_jx$ is also an $SL_2$-invariant.
Hence it is sufficient to show that all products
$w_{j_1}\cdots w_{j_p}$ of degree $2(n+1)$ and all
$x\omega_jy-y\omega_jx$ are linearly independent.
\par
We introduce the lexicographic ordering on $K\langle x,y\rangle$ assuming that
$x<y$. Then by induction we prove that
the minimal monomials $z_{k_1}\cdots z_{k_{2n+2}}$, $z_k\in\{x,y\}$,
of $w_{j_1}\cdots w_{j_p}$ and
$x\omega_jy-y\omega_jx$ have the property that the number of $x$'s in
every beginning $z_{k_1}\cdots z_{k_q}$ of $z_{k_1}\cdots z_{k_{2n+2}}$
is bigger or equal to the number of $y$'s. For example, the minimal monomial of
$w_2=x[x,y]y-y[x,y]x$ is $xxyy$, all its beginnings are
$x,xx,xxy,xxyy$ and the number of entries of $x$ and $y$ are
$(1,0),(2,0),(2,1),(2,2)$, respectively. Similarly, the minimal monomial of
\[
w_1w_2^2=[x,y](x[x,y]y-y[x,y]x)(x[x,y]y-y[x,y]x)
\]
is $xyxxyyxxyy$ and the entries of $x$ and $y$ in the beginnings are
\[
(1,0),(1,1),(2,1),(3,1),(3,2),(3,3),(4,3),(5,3),(5,4),(5,5).
\]
Pay attention that the first place where the number of $x$'s is equal to
the number of $y$'s, namely the beginning $xy$, corresponds to the beginning
$w_1=[x,y]$ in $w_1w_2^2$ and the rest of the minimal monomial
$xxyyxxyy$ has the same property.
\par
We shall show that the products
$w_{j_1}\cdots w_{j_p}$ (including the case $p=1$ of a product of one free generator
$x\omega_jy-y\omega_jx$) are in a 1-1 correspondence with
the words $z_{k_1}\cdots z_{k_{2n+2}}$
in $x$ and $y$ with the property that the number of $x$'s in
every beginning $z_{k_1}\cdots z_{k_q}$
is bigger or equal to the number of $y$'s.
Let $\omega=w_{j_1}\cdots w_{j_p}$ be a product of elements of the constructed set.
If $p=1$, i.e. $\omega=w_j$ is in the set, then $w_j=x\omega'y-y\omega'x$
and the minimal monomial $z_1\cdots z_{2n}$ of $\omega'$ has the property
that the number of $x$'s in every beginning of $z_1\cdots z_{2n}$
is bigger or equal to the corresponding number of $y$'s. Since the
minimal monomial of $w_j$ is $xz_1\cdots z_{2n}y$, we obtain that in every
of its proper beginnings the number of occurances of $x$ is strictly bigger than
the number of entries of $y$. If $p>1$, then,
reading the minimal word from left to right, the first place where
the numbers of the $x$'s and the $y$'s is the same, is the end of $w_{j_1}$.
This arguments combined with induction
easily imply that the different products $w_{j_1}\cdots w_{j_p}$
have different minimal monomials and each word corresponds to some product
$w_{j_1}\cdots w_{j_p}$. Hence the products $w_{j_1}\cdots w_{j_p}$ are linearly
independent and this completes the proof.
\end{proof}

\begin{corollary} \label{constants of Engel algebras}
For any variety $\mathfrak W$ of associative algebras which does not contain the
metabelian variety $\mathfrak M$, the algebra of constants
$F_2({\mathfrak W})^{\delta}$ is finitely generated.
\end{corollary}

\begin{proof}
It is well known that any variety $\mathfrak W$ which does not contain
$\mathfrak M$ satisfies some Engel identity $[x_2,x_1,\ldots,x_1]=0$.
By a theorem of Latyshev \cite{La}
any finitely generated PI-algebra
satisfying a non-matrix polynomial identity, satisfies also
some identity of the form $[x_1,x_2]\cdots[x_{2k-1},x_{2k}]=0$.
Applying this result to $F_2({\mathfrak W})$ we obtain that
$F_2({\mathfrak W})$ is solvable as a Lie algebra,
and, by a theorem of Higgins \cite{Hi}
$F_2({\mathfrak W})$ is Lie nilpotent.
(Actually Zelmanov \cite{Z} proved the stronger result that any
Lie algebra over a field of characteristic zero
satisfying the Engel identity is nilpotent.)
\par
By Drensky \cite{D2}, for any nilpotent variety $\mathfrak W$,
and for a fixed positive integer $m$,
the vector space $B_m({\mathfrak W})$
of so called proper polynomials in $F_m({\mathfrak W})$
is finite dimensional. Using the relation
\[
F_m({\mathfrak W})\cong K[x_1,\ldots,x_m]\otimes_K B_m({\mathfrak W})
\]
between the $GL_m$-module structure of $F_m({\mathfrak W})$
and $B_m({\mathfrak W})$ and the Young rule, we can derive the following.
There exists a positive constant $p$ such that the
nonzero irreducible components $W(\lambda_1,\ldots,\lambda_m)$
of the $GL_m$-module $F_m({\mathfrak W})$ satisfy the restriction
$\lambda_2\leq p$. Hence
the subalgebra $F_2({\mathfrak W})^{SL_2}$ of $SL_2$-invariants of
$F_2({\mathfrak W})$ (which is spanned on the highest weight vectors
of $W(\lambda_1,\lambda_2)$ with $\lambda_1\leq p$) is finite dimensional.
Now the statement follows from Corollary
\ref{generatation by sl-invariants for algebras of rank two}
because $F_2({\mathfrak W})^{\delta}$ is generated by $x$ and
the finite dimensional vector space $F_2({\mathfrak W})^{SL_2}$.
\end{proof}

Corollary \ref{constants of Engel algebras} inspires the following:

\begin{question}
Is it true that, for $m\geq 2$ and for
a fixed nonzero Weitzenb\"ock derivation $\delta$,
the algebra of constants $F_m({\mathfrak W})^{\delta}$
is finitely generated if and only if the variety of
associative algebras $\mathfrak W$ does not contain
the metabelian variety $\mathfrak M$?
\end{question}

Corollary \ref{lifting of constants},
Example \ref{free metabelian associative algebra} and
Corollary \ref{constants of Engel algebras} show that the answer to
this question is affirmative for $m=2$. In the next section we
shall show that the algebra of constants $F_m({\mathfrak W})^{\delta}$
is not finitely generated if $\mathfrak W$ contains $\mathfrak M$.

\section{Constants of Relatively Free Associative Algebras}

First we shall work in the free metabelian associative algebra
$F_m({\mathfrak M})$ where the metabelian variety is defined by
the polynomial identity $[x_1,x_2][x_3,x_4]=0$. We need an embedding of
$F_m({\mathfrak M})$ into a wreath product. For this purpose,
let $Y=\{y_1,\ldots,y_m\}$, $U=\{u_1,\ldots,u_m\}$ and $V=\{v_1,\ldots,v_m\}$
be three sets of commuting variables and let
\[
M=\sum_{i=1}^ma_iK[U,V]
\]
be the free $K[U,V]$-module of rank $m$ generated by $\{a_1,\ldots,a_m\}$.
Clearly, $M$ has also a structure of a free $K[Y]$-bimodule with the action
of $K[Y]$ defined by
\[
y_ja_i=a_iu_j,\quad a_iy_j=a_iv_j,\quad i,j=1,\ldots,m.
\]
Define the trivial multiplication $M\cdot M=0$ on $M$ and
consider the algebra
\[
W=K[Y]\rightthreetimes M,
\]
which is similar to the abelian wreath product of Lie algebras,
see \cite{Sh2} ($M$ is an ideal of $W$ with multiplication by $K[Y]$
induced by the bimodule action of $K[Y]$ on $M$).
Obviously $W$ satisfies the metabelian identity and hence belongs to
$\mathfrak M$.
The following proposition is a partial case of the main result of
Lewin \cite{Le}, see also Umirbaev \cite{U} for further applications of this
construction to automorphisms of relatively free associative algebras.

\begin{proposition} \label{embedding in wreath products}
The mapping $\iota:x_j\to y_j+a_j$, $j=1,\ldots,m$, defines an embedding
$\iota$ of $F_m({\mathfrak M})$ into $W=K[Y]\rightthreetimes M$.
\end{proposition}

\begin{proposition} \label{metabelian algebras of any rank}
For any nontrivial Weitzenb\"ock derivation $\delta$ of
the free metabelian associative algebra $F_m({\mathfrak M})$ of rank
$m\geq 2$, the algebra of constants $F_m({\mathfrak M})^{\delta}$
is not finitely generated.
\end{proposition}

\begin{proof} The derivation $\delta$ acts as a linear operator on the
vector space with basis $\{x_1,\ldots,x_m\}$ and we define in a similar way
the action of $\delta$ on the vector spaces with bases
$\{y_1,\ldots,y_m\}$ and $\{a_1,\ldots,a_m\}$:
If $\delta(x_j)=\sum_{i=1}^m\alpha_{ij}x_j$, $\alpha_{ij}\in K$, then
$\delta(y_j)=\sum_{i=1}^m\alpha_{ij}y_j$ and
$\delta(a_j)=\sum_{i=1}^m\alpha_{ij}a_j$, $j=1,\ldots,m$.
As in the proof of Lemma \ref{extending derivations} we can show that
this action $\delta$ defines a derivation on $W$ and on the polynomial
algebra $K[U,V]$ (which we denote also by $\delta$).
Additionally, we consider the embedding $\iota$
of $F_m({\mathfrak M})^{\delta}$
as a subalgebra in $W$, as stated in
Proposition \ref{embedding in wreath products}. By definition
$\delta(\iota(x_j))=\delta(y_j+a_j)=\iota(\delta(x_j))$ and hence if
$\delta(f(X))=0$ in $F_m({\mathfrak M})$, then
the same holds for the image $\iota(f)$ of $f$ in $W$.
In this way, $\iota$ embeds the algebra
of constants $F_m({\mathfrak M})^{\delta}$ into
the algebra of constants $W^{\delta}$.
\par
As till now, we assume that $\delta(x_1)=0$ and $\delta(x_2)=x_1$. If the
algebra of constants $F_m({\mathfrak M})^{\delta}$ is generated
by a finite set $\{f_1,\ldots,f_n\}$, then, as elements of $W$,
\[
\iota(f_k)=g_k(Y)+\sum_{i=1}^ma_ih_{ik}(U,V),\quad
g_k(Y)\in K[Y], h_{ik}(U,V)\in K[U,V],
\]
$i=1,\ldots,m$, $k=1,\ldots,n$, and
\[
g_k(Y),b_k=\sum_{i=1}^ma_ih_{ik}(U,V),\quad k=1,\ldots,n,
\]
are also constants.
Hence $\iota\left(F_m({\mathfrak M})^{\delta}\right)$
is a subalgebra of the subalgebra of $W^{\delta}$
generated by the union of the finite sets
\[
\{g_1,\ldots,g_n\}\subset K[Y]^{\delta},\quad
\{b_1,\ldots,b_n\}\subset M^{\delta}.
\]
This implies that $F_m({\mathfrak M})^{\delta}$ is a subalgebra of
\[
W_0=K[Y]^{\delta}\rightthreetimes \sum_{k=1}^nb_kK[U]^{\delta}K[V]^{\delta}.
\]
By Corollary \ref{transcendence degree} the transcendence degree of
$K[Y]^{\delta}$ is equal to $m-1$ and hence the transcendence
degree of $K[U]^{\delta}K[V]^{\delta}$ is equal to $2(m-1)$.
Since, see e.g. the book by
Krause and Lenagan \cite{KL}, the Gelfand-Kirillov dimension of
a commutative algebra is equal to the transcendence degree of the algebra,
we easily derive that the Gelfand-Kirillov dimension of the
algebra $W_0$ is bounded from above by $2(m-1)$.
On the other hand, the vector space $\iota\left([x_1,x_2]\right)K[U,V]$
is contained in $\iota\left(F_m({\mathfrak M})\right)$ and is
a free $K[U,V]$-module generated by $a_1(v_2-u_2)+a_2(u_1-v_1)$.
Since $\iota\left([x_1,x_2]\right)\in M^{\delta}$, we obtain that
$\iota\left([x_1,x_2]\right)K[U,V]^{\delta}$ is a free
$K[U,V]^{\delta}$-module. By Corollary \ref{transcendence degree}
the transcendence degree of $K[U,V]^{\delta}$ is equal to $2m-1$,
and hence the Gelfand-Kirillov dimension of the $K[U,V]^{\delta}$-module
is equal to $2m-1$. This is also a lower bound for the Gelfand-Kirillov dimension
of $F_m({\mathfrak M})^{\delta}$ which contradicts with the inequality
$\text{\rm GKdim}(F_m({\mathfrak M})^{\delta})\leq \text{\rm GKdim}(W_0)\leq 2(m-1)$.
\end{proof}

\begin{remark} \label{concrete elements for metabelian case}
In the notation of Proposition \ref{metabelian algebras of any rank},
if $b_1,\ldots,b_k$ is a finite number of elements in $M^{\delta}$, then
the subalgebra of $K[Y]^{\delta}\rightthreetimes M^{\delta}$
generated by $K[Y]^{\delta}$ and $b_1,\ldots,b_k$, contains only
a finite number of elements $\iota([x_1,x_2])(u_1v_2-u_2v_1)^n$. This can be seen
in the following way. We consider
the localization of the polynomial algebra
$K[Y][y_1^{-1}]=K[y_1,y_2,\ldots,y_m][y_1^{-1}]$, and similarly
$K[U][u_1^{-1}],K[V][v_1^{-1}]$. Then we define
$W'=K[Y][y_1^{-1}]\rightthreetimes MK[u_1^{-1},v_1^{-1}]$. Since
$y_1,u_1,v_1$ are $\delta$-constants, we can extend the action of
$\delta$ as a derivation on $W$ to a derivation on $W'$.
Let $\delta^{p_j+1}(y_j)=0$, $j=1,\ldots,m$, and let us define
\[
\tilde y_j=\sum_{k=0}^{p_j}\frac{\delta^k(y_j)}{k!}(-y_2)^ky_1^{p_j-k},\quad
j=3,4,\ldots,m,
\]
and similarly $\tilde y_j, \tilde u_j, \tilde v_j$. Let also
$\tilde w_2=u_1v_2-u_2v_1$. By Proposition \ref{generators Weitzenboeck}
\[
(K[Y][y_1^{-1}])^{\delta}=K[y_1,y_1^{-1}][\tilde y_3,\tilde y_4,\ldots,\tilde y_m],
\]
\[
(K[U,V][u_1^{-1},v_1^{-1}])^{\delta}=K[u_1,v_1,u_1^{-1},v_1^{-1}]
[\tilde u_3,\ldots,\tilde u_m,\tilde v_3,\ldots,\tilde u_m,\tilde w_2].
\]
The algebra generated by $K[Y]^{\delta}$ and $b_1,\ldots,b_k$
is a subalgebra of
\[
(K[Y][y_1^{-1}])^{\delta}\rightthreetimes
\sum_{j=1}^kb_j(K[U][u_1^{-1}])^{\delta}(K[V][v_1^{-1}])^{\delta}
\]
and hence its elements have the form
\[
f(\tilde y_3,\ldots,\tilde y_m)
+\sum_{j=1}^mb_jf_j(\tilde u_3,\ldots,\tilde u_m,\tilde v_3,\ldots,\tilde v_m),
\]
where $f$ and $f_j$ are polynomials with coefficients depending
respectively on $y_1,y_1^{-1}$ and $u_1,v_1,u_1^{-1},v_1^{-1}$.
Since $\tilde u_3,\ldots,\tilde u_m,\tilde v_3,\ldots,\tilde u_m,\tilde w_2$
are algebraically independent on $K[u_1,v_1,u_1^{-1},v_1^{-1}]$,
and the finite number of elements $b_1,\ldots,b_k$
contains only a finite number of summands, we cannot present all elements
$\iota([x_1,x_2])(u_1v_2-u_2v_1)^n=(a_1(v_2-u_2)+a_2(u_1-v_1))\tilde w_2^n$
in the form
\[
(a_1(v_2-u_2)+a_2(u_1-v_1))\tilde w_2^n=
\sum_{j=1}^mb_jf_{jn}(\tilde u_3,\ldots,\tilde u_m,\tilde v_3,\ldots,\tilde v_m).
\]
\end{remark}

\begin{theorem} \label{constants of relatively free associative algebras}
Let $\mathfrak W$ be a variety of associative algebras containing
the metabelian variety $\mathfrak M$. Then for any
$m\geq 2$ and for any fixed nonzero Weitzenb\"ock derivation $\delta$,
the algebra of constants $F_m({\mathfrak W})^{\delta}$
is not finitely generated.
\end{theorem}

\begin{proof}
By Corollary \ref{lifting of constants}
the algebra $F_m({\mathfrak M})^{\delta}$
is a homomorphic image of $F_m({\mathfrak W})^{\delta}$.
Now the proof follows immediately because
$F_m({\mathfrak M})^{\delta}$
is not finitely generated by
Proposition \ref{metabelian algebras of any rank}.
\end{proof}

\begin{remark}
Using the elements $\iota([x_1,x_2])(u_1v_2-u_2v_1)^n$, $n\geq 0$, from
Remark \ref{concrete elements for metabelian case}
for any variety $\mathfrak W$ containing the metabelian variety
$\mathfrak M$ and any nontrivial Weitzenb\"ock derivation $\delta$
we can construct an infinite set of constants which is not contained
in any finitely generated subalgebra of $F_m({\mathfrak W})^{\delta}$.
Again, we assume that
$\delta(x_1)=0$, $\delta(x_2)=x_1$.
Let $l_u$ and $r_u$ be, respectively, the operators of left and right
multiplication by $u\in F_m({\mathfrak W})$. Consider the elements
\[
(l_{x_1}r_{x_2}-l_{x_2}r_{x_1})^n[x_1,x_2],\quad
n\geq 0.
\]
All these elements are constants which are liftings of
the constants from Remark \ref{concrete elements for metabelian case} and hence
any finitely generated subalgebra of $F_m({\mathfrak W})^{\delta}$
does not contain $(l_{x_1}r_{x_2}-l_{x_2}r_{x_1})^n[x_1,x_2]$ for
sufficiently large $n$.
\end{remark}

\begin{corollary} \label{unitriangular invariants}
Let $\mathfrak W$ be a variety of associative algebras containing
the metabelian variaty $\mathfrak M$. Then for any $m\geq 2$
the algebra $F_m({\mathfrak W})^{UT_m}$ of $UT_m(K)$-invariants
is not finitely generated.
\end{corollary}

\begin{proof}
Let the algebra $F_m({\mathfrak W})^{UT_m}$ be finitely generated.
By Remarks \ref{three remarks}, the algebra
$(K\langle x_1,\ldots,x_m\rangle)^{UT_m}$, and hence
also $F_m({\mathfrak W})^{UT_m}$
is spanned by all highest weight vectors.
Hence $F_m({\mathfrak W})^{UT_m}$ is generated by a finite system
of highest weight vectors $w(x_1,\ldots,x_m)\in W(\lambda)\subset
F_m({\mathfrak W})^{UT_m}$.
Hence $F_m({\mathfrak W})^{UT_m}$ is multigraded and has
a finite multihomogeneous set of generators. The generators
which depend on $x_1$ and $x_2$ only, generate the subalgebra
spanned by all highest weight vectors
$w(x_1,\ldots,x_m)\in W(\lambda_1,\lambda_2,0,\ldots,0)$.
This subalgebra coincides with the algebra of $UT_2$-invariants
of $F_2({\mathfrak W})$ and hence with the algebra of constants
of the Weitzenb\"ock derivation $\delta$ of $F_2({\mathfrak W})$
defined by $\delta(x_1)=0$, $\delta(x_2)=x_1$. By
Theorem \ref{constants of relatively free associative algebras} for $m=2$
(or by Corollary \ref{lifting of constants} and
Example \ref{free metabelian associative algebra})
$F_2({\mathfrak W})^{\delta}$ is not finitely generated.
Hence the algebra $F_m({\mathfrak W})^{UT_m}$ cannot be finitely generated.
\end{proof}

\begin{remark} \label{unitriangular invariants of Lie nilpotent varieties}
Let $\mathfrak W$ be a Lie nilpotent variety of associative algebras and let
$m$ be a fixed positive integer.
Using the approach of \cite{D2} (as in the proof of
Corollary \ref{constants of Engel algebras}),
and the fact that $F_m({\mathfrak W})$ is a direct sum of
$GL_m$-modules of the form $W(\lambda_1,\ldots,\lambda_m)$ with
$\lambda_2\leq p$ for some $p$, one can
show that there exists a finite system of highest weight
vectors $w_i(x_1,\ldots,x_k)\in F_m({\mathfrak W})$, $i=1,\ldots,k$,
such that all highest weight vectors of $F_m({\mathfrak W})$
are linear combinations of $x^nw_i(x_1,\ldots,x_k)$. Hence
the algebra $F_m({\mathfrak W})^{UT_m}$ of $UT_m$-invariants
is generated by $x$ and $w_i(x_1,\ldots,x_k)$, $i=1,\ldots,k$.
Hence $F_m({\mathfrak W})^{UT_m}$ is finitely generated.
\end{remark}

\section{Generic $2\times 2$ Matrices}

In this section we construct classes
of automorphisms of the relatively free algebra
$F_2(\text{\rm var }M_2(K))$. This algebra is isomorphic to the
algebra generated by two generic $2\times 2$ matrices $x$ and $y$.
So, the results are stated in the natural setup of the trace algebra.
We start with the necessary background, see Formanek \cite{F2},
Alev and Le Bruyn \cite{AL}, or Drensky and Gupta \cite{DG}.
\par
We consider the polynomial algebra in 8 variables
$\Omega=K[x_{ij},y_{ij}\mid i,j=1,2]$.
The algebra $R$ of two generic $2\times 2$ matrices
\[
x = \begin{pmatrix}
x_{11}&x_{12}\\
x_{21}&x_{22}\\
\end{pmatrix}\quad
\text{ \rm and }\quad
y = \begin{pmatrix}
y_{11}&y_{12}\\
y_{21}&y_{22}\\
\end{pmatrix}
\]
is the subalgebra of $M_2(\Omega)$ generated by $x$ and $y$. We denote
by $C$ the centre of $R$ and by $\bar C$ the
algebra generated by all the traces of elements from $R$. Identifying
the elements of $\bar C$ with $2 \times 2$ scalar matrices we
denote by $T$ the generic trace algebra generated by $R$ and $\bar C$.
It is well known that $\bar C$ is generated by
\[
\text{\rm tr}(x), \text{\rm tr}(y), \text{\rm det}(x), \text{\rm det}(y),
\text{\rm tr}(xy)
\]
and is isomorphic to the polynomial algebra
in five variables.

\begin{proposition} \label{centre of generic matrix algebra}
{\rm (Formanek, Halpin, Li \cite{FHL})}
The vector subspace of $C$
consisting of all polynomials without constant term is a free
$\bar C$-module generated by $[x,y]^2$.
\end{proposition}

For our purposes it is more convenient to replace in $T$
(as in \cite{AL}) the generic
matrices $x$ and $y$ by the generic traceless matrices
$$
x_0=x-\frac{1}{2}\text{\rm tr}(x),\,
y_0=y-\frac{1}{2}\text{\rm tr}(y)
$$
and assume that $T$ is generated by $x_0$, $y_0$, $\text{\rm tr}(x)$,
$\text{\rm tr}(y)$, $\text{\rm det}(x_0)$,
$\text{\rm det}(y_0)$, $\text{\rm tr}(x_0y_0)$. A further reduction is
to use the formulas
\[
\text{\rm det}(x_0)=
-\frac{1}{2}\text{\rm tr}(x_0^2),\,
\text{\rm det}(y_0)=
-\frac{1}{2}\text{\rm tr}(y_0^2),
\]
and to replace the determinants by $\text{\rm tr}(x_0^2)$ and
$\text{\rm tr}(y_0^2)$.
In this way, we may assume that $\bar C$ is
generated by
\[
p = \text{\rm tr}(x), q = \text{\rm tr}(y), u = \text{\rm tr}(x_0^2),
v = \text{\rm tr}(y_0^2), t = \text{\rm tr}(x_0y_0).
\]
Then $[x,y]^2=t^2-uv$ and
\[
T=\bar C+\bar Cx_0+\bar Cy_0+\bar C[x_0,y_0]
\]
is a free $\bar C$-module generated by $1,x_0,y_0,[x_0,y_0]$.
\par
The defining relations of the algebra generated by the $2\times 2$ traceless
matrices $x_0$ and $y_0$ are $[x_0^2,y_0]=[y_0^2,x_0]=0$,
see e.g. \cite{LB} or \cite{DKo} for the case of characteristic 0
and \cite{Ko2} for the case of an arbitrary infinite base field.
More generally, the defining relations of the algebra generated by
$m$ generic $2\times 2$ traceless matrices $y_1,\ldots,y_m$
are $[v_1^2,v_2]=0$, where $v_1,v_2$ run on the set of all Lie elements in
$K\langle y_1,\ldots,y_m\rangle$ which is a restatement
of the theorem of Razmyslov \cite{R} for the weak polynomial
identities of $M_2(K)$.
An explicitly written system of defining relations consists of
$[y_i^2,y_j]=0$, $i,j=1,\ldots,m$, and the standard polynomials
$s_4(y_{i_1},y_{i_2},y_{i_3},y_{i_4})=0$, $1\leq i_1<i_2<i_3<i_4\leq m$,
see \cite{DKo}.

\begin{lemma} \label{derivations of generic matrices}
Every mapping $\delta: \{p,q,x_0,y_0\}\to T$ such that
\[
\delta(p),\delta(q)\in \bar C,\quad
\delta(x_0),\delta(y_0)\in \bar Cx_0+\bar Cy_0+\bar C[x_0,y_0]
\]
can be extended to a derivation of $T$.
\end{lemma}

\begin{proof} The defining relations of $T$ are
\[
[p,q]=[p,x_0]=[p,y_0]=[q,x_0]=[q,y_0]=0,
\]
together with the defining relations of the subalgebra generated by $x_0,y_0$.
It is sufficient to see that the extension of $\delta$
(inductively, by the rule $\delta(fg)=\delta(f)g+f\delta(g)$)
to a derivation on $T$
is well defined, i.e. sends the defining relations to 0. For the relations involving
$p$ and $q$ this can be checked directly:
\[
\delta([p,q])=[\delta(p),q]+[p,\delta(q)]=0,
\]
analoguously for $\delta([p,x_0]),\delta([p,y_0]),\delta([q,x_0]),\delta([q,y_0])$,
because $p,q,\delta(p),\delta(q)$ are in the centre of $T$.
The condition for the defining relations of the algebra generated by $x_0,y_0$
can be proved using the universal properties of this algebra or directly:
Since $x_0^2, y_0^2, x_0y_0+y_0x_0,[x_0,y_0]^2$ are in the centre of $T$,
and $x_0[x_0,y_0]+[x_0,y_0]x_0=y_0[x_0,y_0]+[x_0,y_0]y_0=0$,
if $\delta(x_0)=ax_0+by_0+c[x_0,y_0]$, $a,b,c\in\bar C$, then
\[
(\delta(x_0))^2=a^2x_0^2+b^2y_0^2+c^2[x_0,y_0]^2+ab(x_0y_0+y_0x_0),
\]
\[
\delta(x_0)x_0+x_0\delta(x_0)=ax_0^2+b(x_0y_0+y_0x_0)
\]
are in the centre of $T$ and $\delta([x_0^2,y_0])=0$. In the same way
$\delta([y_0^2,x_0])=0$.
\end{proof}

\begin{example}\label{automorphisms fixing x}
Let us consider the basic Weitzenb\"ock derivation $\delta$
defined on the relatively free algebra
$F_2(\text{\rm var }M_2(K))$ in its realization as the generic trace algebra
generated by generic $2\times 2$ matrices $x$ and $y$
by $\delta(x)=0$, $\delta(y)=x$. We extend $\delta$ to the trace algebra $T$
by
\[
\delta(p)=\delta(\text{\rm tr}(x))=\text{\rm tr}(\delta(x)),
\]
\[
\delta(q)=\delta(\text{\rm tr}(y))=\text{\rm tr}(\delta(y)),
\]
\[
\delta(x_0)=0,\quad \delta(y_0)=x_0,
\]
\[
\delta(u)=\delta(\text{\rm tr}(x_0^2))=\text{\rm tr}(\delta(x_0^2)),
\]
\[
\delta(v)=\delta(\text{\rm tr}(y_0^2))=\text{\rm tr}(\delta(y_0^2))
\]
\[
\delta(t)=\delta(\text{\rm tr}(x_0y_0))=\text{\rm tr}(\delta(x_0y_0)).
\]
By Lemma  \ref{derivations of generic matrices} this is possible.
Direct calculations give that
\[
\delta(p)=0,\quad \delta(q)=p,\quad
\delta(u)=0,\quad \delta(t)=u,\quad \delta(v)=2t.
\]
Replacing $v$ with $2v_1$, we obtain that the action of
$\delta$ on $\bar C=K[p,q,u,t,v_1]$ is
as in Examples \ref{examples of concrete generators}. Hence
\[
(\bar C)^{\delta}=K[p,u,pt-qu,t^2-2uv_1,2p^2v_1-2pqt+q^2u]
\]
\[
=K[p,u,pt-qu,t^2-uv,q^2u-2pqt+p^2v].
\]
The generators of $(\bar C)^{\delta}$ satisfy the relation
\[
u(q^2u-2pqt+p^2v)+p^2(t^2-uv)=(pt-qu)^2.
\]
If $w\in(\bar C)^{\delta}$, then $\exp(w\delta)$ is an automorphism of $T$.
If $t^2-uv$ divides $w$, then $\exp(w\delta)$ is an automorphism also of $R$.
This automorphism acts on $R$ as
\[
\exp(w\delta): x\to x,\quad
\exp(w\delta): y\to y+wx,
\]
where
$w=(t^2-uv)w_1(p,u,pt-qu,t^2-uv,q^2u-2pqt+p^2v)$ for some polynomial $w_1$.
Such automorphisms (fixing $x$) were studied in the Ph. D. Thesis of Chang \cite{C}.
\end{example}

\begin{example}
Now we shall modify Example \ref{automorphisms fixing x} in the following way.
We use Lemma  \ref{derivations of generic matrices} and define the derivation
$\delta$ of $T$ by
\[
\delta(p)=\alpha_1u+\beta_1t+\gamma_1v,
\quad
\delta(q)=p+\alpha_2u+\beta_2t+\gamma_2v,
\]
$\alpha_i,\beta_i,\gamma_i\in \bar C$, $i=1,2$,
\[
\delta(u)=0,\quad \delta(t)=u,\quad \delta(v)=2t.
\]
This derivation is locally nilpotent and acts on the generic matrices
$x=\frac{1}{2}\text{\rm tr}(x)+x_0$ and $y=\frac{1}{2}\text{\rm tr}(y)+y_0$ by
\[
\delta(x)=\frac{1}{2}(\alpha_1u+\beta_1t+\gamma_1v),\quad
\delta(y)=x+\frac{1}{2}(\alpha_2u+\beta_2t+\gamma_2v).
\]
The matrix of the linear operator $\delta$ acting on the vector space
$Kp+Kq+Ku+Kt+Kv$ (with respect to the basis $\{p,q,u,t,v\}$) is
\[
\begin{pmatrix}
0&1&0&0&0\\
0&0&0&0&0\\
\alpha_1&\alpha_2&0&1&0\\
\beta_1&\beta_2&0&0&2\\
\gamma_1&\gamma_2&0&0&0
\end{pmatrix}
\]
and has rank 3 or 4 depending on whether $\gamma_1=0$ or $\gamma_1\not=0$.
Hence its Jordan normal form is one of the following matrices:
\[
\begin{pmatrix}
0&1&0&0&\\
0&0&1&0&\\
0&0&0&1&\\
0&0&0&0&\\
&&&&0
\end{pmatrix},
\quad
\begin{pmatrix}
0&1&0&&\\
0&0&1&&\\
0&0&0&&\\
&&&0&1\\
&&&0&0
\end{pmatrix},
\quad
\begin{pmatrix}
0&1&0&0&0\\
0&0&1&0&0\\
0&0&0&1&0\\
0&0&0&0&1\\
0&0&0&0&0
\end{pmatrix}.
\]
Examples \ref{examples of concrete generators}
give concrete systems of generators of the algebras of constants
of $(\bar C)^{\delta}$ and hence automorphisms of the algebras $T$ and $R$.
\par
For example, if we fix
$\delta(p)=v$, $\delta(q)=p$, then $\delta$ is a basic derivation with
\[
\delta(q)=p,\quad \delta(p)=v,\quad ,\delta(v)=2t,\quad \delta(t)=u,\quad
\delta(u)=0.
\]
Considering $\bar C=K[q/2,p/2,v/2,t,u]$,
we obtain after some easy calculations that
the algebra of constants is generated by
\[
u,\quad t^2-uv,\quad tp-qu-\frac{v^2}{4},
\]
\[
t^3-\frac{3}{2}utv+\frac{3}{2}u^2p,
\quad
3t^2q-\frac{3}{2}tvp+\frac{v^3}{4}-3uvq+\frac{9}{4}up^2.
\]
In this case $\delta$ acts on $x$ and $y$ by
\[
\delta(x)=\frac{1}{2}\text{\rm tr}(y_0^2)=\frac{1}{2}v,\quad
\delta(y)=x.
\]
If $w$ is in $(\bar C)^{\delta}$, then
\[
\exp(w\delta):x\to x+\frac{wv}{2.1!}+\frac{w^2t}{2!}+\frac{w^3u}{3!},
\]
\[
\exp(w\delta):y\to y+\frac{wx}{1!}+\frac{w^2v}{2.2!}+\frac{w^3t}{3!}+\frac{w^4u}{4!}.
\]
If $w$ is divisible by $t^2-uv$, then $\exp(w\delta)$ is also an automorphism of $R$.
Since all these automorphisms $\exp(w\delta)$
are obtained by the construction of Martha Smith \cite{Sm}, they
induce stably tame automorphisms of $\bar C=K[p,q,u,t,v]$.
\end{example}

\section{Relatively Free Lie Algebras}

We start with few examples for the algebras of constants of relatively free algebras.
By the well known dichotomy a variety of Lie algebras either
satisfies the Engel condition (and by the theorem of Zelmanov \cite{Z}
is nilpotent) or contains the metabelian
variety $\mathfrak A^2$ (which consists of all solvable of class 2 Lie algebras
and is defined by the identity $[[x_1,x_2],[x_3,x_4]]=0$). Since the finitely generated
nilpotent Lie algebras are finite dimensional, the problem for the finite generation
of the algebras of constants of relatively free nilpotent Lie algebras is solved trivially.
\par
The bases of the free polynilpotent Lie algebras were described by Shmelkin \cite{Sh1}.
Considering relatively free algebras of rank 2, we assume that the algebra is generated
by $x$ and $y$ and the basic Weitzenb\"ock derivation $\delta$
is defined by $\delta(x)=0$, $\delta(y)=x$.

\begin{example}
Let $L_2({\mathfrak A}^2)=L_2/L''_2$ be the free metabelian Lie algebra of rank 2.
It has a basis
\[
\{x,y,[y,x,\underbrace{x,\ldots,x}_{a\text{ \rm times}},
\underbrace{y,\ldots,y}_{b\text{ \rm times}}]\mid a,b\geq 0\}.
\]
It is well known (and can be also obtained by simple argumensts from
the Hilbert series of $L_m({\mathfrak A}^2)$)
that the $n$-th cocharacter of the variety ${\mathfrak A}^2$ is
\[
\chi_1({\mathfrak A}^2)=\chi_{(1)},\quad
\chi_n({\mathfrak A}^2)=\chi_{(n-1,1)},\, n\geq 2.
\]
The corresponding highest weight vectors are
\[
w_{(1)}=x,\quad
w_{(n-1,1)}=[y,x,\underbrace{x,\ldots,x}_{n-2\text{ \rm times}}],\, n\geq 2.
\]
Hence the algebra of constants  $L_2({\mathfrak A}^2)^{\delta}$ is generated by
$x$ and $[x,y]$.
\end{example}

\begin{example}
The free abelian-by-\{nilpotent of class 2\} Lie algebra
$L_2({\mathfrak A}{\mathfrak N}_2)=L_2/[L_2,L_2,L_2]'$
satisfies the identity
\[
[[x_1,x_2,x_3],[x_4,x_5,x_6]]=0
\]
and has a basis
\[
\{x,y,[x,y],[y,x,\underbrace{x,\ldots,x}_{a\text{ \rm times}},
\underbrace{y,\ldots,y}_{b\text{ \rm times}},
\underbrace{[x,y],\ldots,[x,y]}_{c\text{ \rm times}}]\mid a+b>0,c\geq 0\}.
\]
Its Hilbert series is
\[
H(L_2({\mathfrak A}{\mathfrak N}_2),t_1,t_2)=t_1+t_2+t_1t_2
+\frac{t_1t_2(t_1+t_2)}{(1-t_1)(1-t_2)(1-t_1t_2)}
\]
\[
=S_{(1)}(t_1,t_2)+S_{(1^2)}(t_1,t_2)
+\sum_{\lambda_1>\lambda_2\geq 1}S_{(\lambda_1,\lambda_2)}(t_1,t_2)
\]
and the highest weight vectors of $L_2({\mathfrak A}{\mathfrak N}_2)$ are
\[
x,\quad [x,y],\quad [y,x,\underbrace{x,\ldots,x}_{a\text{ \rm times}},
\underbrace{[x,y],\ldots,[x,y]}_{c\text{ \rm times}}],\ a>0,c\geq 0.
\]
Hence the algebra $L_2({\mathfrak A}{\mathfrak N}_2)^{\delta}$
is generated by $x$ and $[x,y]$.
\end{example}

\begin{example}
We consider the relatively free algebra $L_2(\text{\rm var }sl_2(K))$
of the variety of Lie algebras generated by the algebra of $2\times 2$
traceless matrices. This algebra is isomorphic to the Lie algebra
generated by the generic $2\times 2$ traceless matrices $x_0,y_0$
considered in Section 7. By Drensky \cite{D1}, as a $GL_2$-module
$L_2(\text{\rm var }sl_2(K))$ has the decomposition
\[
L_2(\text{\rm var }sl_2(K))\cong W(1)\bigoplus\sum W(\lambda_1,\lambda_2),
\]
where the summation runs on all $\lambda=(\lambda_1,\lambda_2)$ such that
$\lambda_2>0$ and at least one of the integers $\lambda_1,\lambda_2$ is odd.
The highest weight vectors of $W(\lambda_1,\lambda_2)$
are given in \cite{D1} but we do not need their concrete form
for our purposes. The algebra of constants $L_2(\text{\rm var }sl_2(K))^{\delta}$
is bigraded. Assuming that the degree of $x$ corresponds to $t$ and the degree of
$y$ is $u=v/t$, the Hilbert series of $L_2(\text{\rm var }sl_2(K))^{\delta}$ is
\[
H(L_2(\text{\rm var }sl_2(K))^{\delta},t,v)=
t+v\left(\sum_{p,q\geq 0}t^pv^q-\sum_{p,q\geq 0}t^{2p}v^{2q+1}\right)
\]
\[
=t+\frac{v}{(1-t)(1-v)}-\frac{v^2}{(1-t)^2(1-v)^2}.
\]
If $L_2(\text{\rm var }sl_2(K))^{\delta}$ is finitely generated, we may
fix a finite system of bigraded generators. For every
homogeneous $f\in L_2(\text{\rm var }sl_2(K))^{\delta}$ we have
$\text{deg}_xf\geq\text{deg}_yf$. Hence the subalgebra
spanned on the homogeneous components of bidegree $(n,n)$, $n$ odd,
is also finitely generated. This subalgebra
is infinite dimensional and its Hilbert series is
obtained from the Hilbert series $H(L_2(\text{\rm var }sl_2(K))^{\delta},t,v)$
by the substitution $t=0$, i.e.
\[
H(L_2(\text{\rm var }sl_2(K))^{\delta},0,v)
=\frac{v}{1-v}-\frac{v^2}{(1-v)^2}.
\]
Besides, the subalgebra is abelian because the commutator of any
two highest weight vectors $w_{(2p+1,2p+1)}$ and
$w_{(2q+1,2q+1)}$ is a highest weight vector
$w_{(2(p+q+1),2(p+q+1))}$ which does not participate in the decomposition of
$L_2(\text{\rm var }sl_2(K))^{\delta}$.
Since the finitely generated abelian Lie algebras are finite dimensional,
we obtain a contradiction which gives that
$L_2(\text{\rm var }sl_2(K))^{\delta}$
cannot be finitely generated.
\end{example}

\begin{example} \label{abelian-by-nilpotent of class 3}
The free abelian-by-\{nilpotent of class 3\} Lie algebra
$L_2({\mathfrak A}{\mathfrak N}_3)=L_2/[L_2,L_2,L_2,L_2]'$ has a basis
consisiting of $x,y$ and commutators of the form
\[
[y,x,\underbrace{x,\ldots,x}_{a\text{ \rm times}},
\underbrace{y,\ldots,y}_{b\text{ \rm times}},
\underbrace{[x,y],\ldots,[x,y]}_{c\text{ \rm times}},
\underbrace{[y,x,x],\ldots,[y,x,x]}_{d\text{ \rm times}},
\underbrace{[y,x,y],\ldots,[y,x,y]}_{e\text{ \rm times}}],
\]
with some natural restrictions of $a,b,c,d,e\geq 0$
which guarantee that these commutators are different from zero
and, up to a sign, pairwise different.
If the algebra of constants $L_2({\mathfrak A}{\mathfrak N}_3)^{\delta}$
is finitely generated, then it has a generating set consisting of
a finite number of bihomogeneous elements $w_1,\ldots,w_k$
of degree $\geq 4$ (and bidegree $(n_1,n_2)$, where $n_1\geq n_2$)
and constants of degree $\leq 3$ (i.e. $x,[x,y],[y,x,x]$).
Since the commutators of length $\geq 4$ commute, we derive that
$L_2({\mathfrak A}{\mathfrak N}_3)^{\delta}$ is a sum of the Lie subalgebra
$N$ generated by $x,[x,y],[y,x,x]$ and the $N$-module generated by
$w_1,\ldots,w_k$. The following elements are constants:
\[
u_n=\sum_{\rho,\sigma,\ldots,\tau\in S_2}\text{\rm sign}(\rho\sigma\cdots\tau)
[y,x,x,z_{\rho(1)},z_{\sigma(1)},\ldots,z_{\tau(1)},
\]
\[
[x,y,z_{\rho(1)}],[x,y,z_{\sigma(1)}],\ldots,[x,y,z_{\tau(1)}]],
\]
where $\{z_1,z_2\}=\{x,y\}$ and, in the summation,
$\rho,\sigma,\ldots,\tau$ run on $n$ copies of
the symmetric group $S_2$. They are homogeneous of bidegree
$(2n+2,2n+1)$ and hence can be written as linear combinations of commutators
involing a $w_i$, several $[x,y]$ and not more than one $x$ or $[y,x,x]$.
But this is impossible because for sufficiently large $n$
one cannot obtain the summands of $u_n$
\[
[y,x,x,\underbrace{x,\ldots,x}_{n\text{ \rm times}},
\underbrace{[x,y,y],\ldots,[x,y,y]}_{n\text{ \rm times}}].
\]
Hence the algebra $L_2({\mathfrak A}{\mathfrak N}_3)^{\delta}$
is not finitely generated.
\end{example}

\begin{example}
Let $m>2$ and let $\delta$ be the Weitzenb\"ock derivation of the free
metabelian Lie algebra $L_m({\mathfrak A}^2)$ defined by
$\delta(x_2)=x_1$, $\delta(x_j)=0$ for $j\not=2$. Then, since
$L_m({\mathfrak A}^2)$ has a basis consisting of $x_j$ and all
commutators $[x_{i_1},x_{i_2},\ldots,x_{i_n}]$
with $i_1>i_2\leq i_3\leq\cdots\leq i_n$, then
the free generators $x_j$, $j\not=2$, and the commutators
which do not include $x_2$ are constants.
It is easy to see that the commutators with $x_2$ are of the form
\[
u'=[x_2,\underbrace{x_1,\ldots,x_1}_{a\text{ \rm times}},
\underbrace{x_2,\ldots,x_2}_{b\text{ \rm times}},x_{i_k},\ldots,x_{i_n}],\quad
a>0,b\geq 0,i_k>2,
\]
\[
u''=[x_i,\underbrace{x_1,\ldots,x_1}_{a\text{ \rm times}},
\underbrace{x_2,\ldots,x_2}_{b\text{ \rm times}},
x_{i_k},\ldots,x_{i_n}],\quad
a\geq 0,b>0,i,i_k>2.
\]
It is easy to see that a linear combination of $u'$ and $u''$ is a constant
if and only if it contains as summands only
$u'$ with $b=0$ and does not contain any $u''$. Hence the algebra of constants
$L_m({\mathfrak A}^2)^{\delta}$ is generated by $x_1,x_j$, $j>2$, and $[x_1,x_2]$.
\end{example}

\begin{example}
Let $m>2$ and let $\delta$ be the Weitzenb\"ock derivation of the free
abelian-by-\{nilpotent of class 2\} Lie algebra
$L_m({\mathfrak A}{\mathfrak N}_2)$ defined, as in the previous example, by
$\delta(x_2)=x_1$, $\delta(x_j)=0$ for $j\not=2$.
We define a $GL_2$-action on $L_m({\mathfrak A}{\mathfrak N}_2)$
assuming that $GL_2$ fixes $x_3,\ldots,x_m$ and acts canonically
on the linear combinations of $x_1,x_2$. Then the subspaces of
$L_m({\mathfrak A}{\mathfrak N}_2)$ which are homogeneous in each variable
$x_3,\ldots,x_m$ are $GL_2$-invariant. This easily implies that
the algebra of constants $L_m({\mathfrak A}{\mathfrak N}_2)^{\delta}$
is multigraded and $\text{\rm deg}_{x_1}f\geq \text{\rm deg}_{x_2}f$
for each multihomogeneous constant $f$. If the algebra
$L_m({\mathfrak A}{\mathfrak N}_2)^{\delta}$ is finitely generated, then
as in Example \ref{abelian-by-nilpotent of class 3}, it is generated
by $x_1,[x_1,x_2],x_3,x_4,\ldots,x_m$ and a finite system
$w_1,\ldots,w_k$ of homogeneous elements of degree $\geq 3$.
Then $L_m({\mathfrak A}{\mathfrak N}_2)^{\delta}$ is a sum of the subalgebra
$N$ generated by $x_1,[x_1,x_2],x_3,x_4,\ldots,x_m$
and the $N$-module generated by $w_1,\ldots,w_k$.
The constants
\[
\sum_{\rho,\sigma,\ldots,\tau\in S_2}\text{\rm sign}(\rho\sigma\cdots\tau)
u_n=[x_1,x_2,x_{\rho(1)},x_{\sigma(1)},\ldots,x_{\tau(1)},
\]
\[
[x_3,x_{\rho(1)}],[x_3,x_{\sigma(1)}],\ldots,[x_3,x_{\tau(1)}]],
\]
where in the summation
$\rho,\sigma,\ldots,\tau$ run on $n$ copies of
the symmetric group $S_2$, are homogeneous of degree $(n+1,n+1,n,0,\ldots,0)$
and arguments as in Example \ref{abelian-by-nilpotent of class 3}
show that this is impossible. Hence the algebra
$L_m({\mathfrak A}{\mathfrak N}_2)^{\delta}$ cannot be finitely generated.
\end{example}

In the above examples, the matrix of the Weitzenb\"ock derivation $\delta$
(as a linear operator acting on the vector space with basis $\{x_1,\ldots,x_m\}$)
is of rank 1. This gives rise to the following natural problem.

\begin{problem} If the matrix of the Weitzenb\"ock derivation $\delta$
is of rank $1$, find the exact frontier where the algebra of constants
$L_m({\mathfrak W})^{\delta}$
becomes finitely generated,
i.e. describe all varieties of Lie algebras $\mathfrak W$ and all integers $m>1$
such that the algebra $L_m({\mathfrak W})^{\delta}$ is finitely generated.
\end{problem}

Finally, we shall give the solution of this problem in the case of rank $\geq 2$.

\begin{theorem}
Let $\mathfrak W$ be a nonnilpotent variety of Lie algebras and let
$\delta$ be a Weitzenb\"ock derivation of the relatively free algebra
$L_m({\mathfrak W})$, $m\geq 3$. If the rank of the matrix of $\delta$
is $\geq 2$, then the algebra of constants $L_m({\mathfrak W})^{\delta}$
is not finitely generated.
\end{theorem}

\begin{proof}
As in the associative case, it is sufficient to establish the theorem for
the metabelian variety of Lie algebras only.
We consider the abelian wreath product of Lie algebras
\[
W_m=(Ky_1\oplus \cdots\oplus Ky_m)\rightthreetimes
\sum_{j=1}^ma_jK[y_1,\ldots,y_m],
\]
where $[y_i,y_j]=[a_if_i,a_jf_j]=0$
and $[a_if_i,y_j]=a_if_iy_j$
($f_i,f_j\in K[y_1,\ldots,y_m]$). Then by the theorem of
Shmelkin \cite{Sh2} the mapping
$\iota:x_j\to a_j+y_j$, $j=1,\ldots,m$, defines an embedding
of the free metabelian Lie algebra $L_m({\mathfrak A}^2)$ into
$W_m$. We assume that $\delta$
is in its normal Jordan form (and $\delta(x_2)=x_1$, $\delta(x_1)=0$).
Hence the fixed part of
$Ky_1\oplus \cdots\oplus Ky_m$ is of dimension $m-\text{\rm rank}(\delta)
\leq m-2$ and is spanned on some free generators
$x_{j_1}=x_1,x_{j_2},\ldots,x_{j_p}$, $p\leq m-2$. If the algebra
$L_m({\mathfrak A}^2)^{\delta}$ is finitely generated, then
it is a sum of $Kx_1\oplus Kx_{j_2}\oplus\cdots\oplus Kx_{j_p}$
and a finitely generated $K[x_1,x_{j_2},\ldots,x_{j_p}]$-submodule
of the commutator ideal $L_m({\mathfrak A}^2)'$. But, as in
the associative case, this is impossible because
the image of this module under $\iota$ should
contain for example $\iota([x_2,x_1])K[y_1,\ldots,y_m]^{\delta}$
and the transcendence degree of $K[y_1,\ldots,y_m]^{\delta}$ is equal to $m-1$.
\par
One can see directly, that if $\delta(x_3)=x_2$, then a finitely generated
subalgebra of $\iota\left(L_m({\mathfrak A}^2)^{\delta}\right)$ cannot contain all
constants
\[
\iota([x_2,x_1])(x_2^2-2x_1x_3)^n,\quad n\geq 0.
\]
Similarly, if $\delta(x_4)=x_3,\delta(x_3)=0$, then
$\iota\left(L_m({\mathfrak A}^2)^{\delta}\right)$ cannot contain all
\[
\iota([x_2,x_1])(x_1x_4-x_2x_3)^n,\quad n\geq 0.
\]
\end{proof}

\section*{Acknowledgements}

This project was carried out when the first author visited the Department
of Mathematics of the University of Manitoba in Winnipeg.
He is very thankful for the kind hospitality and the creative atmosphere.
The first author is also very grateful to Andrzej Nowicki for the useful
discussions on Weitzenb\"ock derivations of polynomial algebras.

\vfill\eject
\end{document}